\documentclass[11pt,english,a4paper]{smfart}

\usepackage[utf8]{inputenc}
\usepackage[T1]{fontenc}
\usepackage[francais]{babel}
\usepackage{lmodern}
\usepackage{smfthm}
\usepackage[headings]{fullpage}
\usepackage{amssymb}
\usepackage{enumerate}

\let\cal\mathcal

\let\hat\widehat
\let\tilde\widetilde
\let\phi\varphi

\let\epsilon\varepsilon

\def\Q{{\bf Q}} 
\def\Z{{\bf Z}}
\def\C{{\bf C}}
\def\N{{\bf N}}

\def\A{{\bf A}}
\def\E{{\bf E}}
\def\B{{\bf B}}
\def\O{{\cal O}}
\def\G{{\cal G}}
\def\F{{\bf F}}

\def\k{{\bf k}}

\def\D{{\bf D}}
\def\M{{\bf M}}

\def\Mat{{\mathrm{Mat}}}

\def\pa{{\mathrm{pa}}}

\def\cbf{{\bf{c}}}
\def\kbf{{\bf{k}}}

\def\Zp{{\Z_p}}
\def\Qp{{\Q_p}}
\def\Cp{{\C_p}}

\def\At{{\tilde{\bf{A}}}}
\def\Atplus{{\tilde{\bf{A}}^+}}
\def\Bt{{\tilde{\bf{B}}}}
\def\Btplus{{\tilde{\bf{B}}^+}}
\def\Et{{\tilde{\bf{E}}}}
\def\Etplus{{\tilde{\bf{E}}^+}}

\def\Gal{{\mathrm{Gal}}}

\def\GL{{\mathrm{GL}}}

\def\k{{\mathbf{k}}}
\def\la{{\mathrm{la}}}
\def\an{{\mathrm{an}}}
\def\cycl{{\mathrm{cycl}}}
\def\rig{{\mathrm{rig}}}

\author{Aditya Karnataki}
\address{BICMR\\
Peking University}
\email{adityack@bicmr.pku.edu.cn}
\urladdr{bicmr.pku.edu.cn/~karnataki/}

\author{Léo Poyeton}
\address{BICMR\\
Peking University}
\email{leo.poyeton@ens-lyon.fr}
\urladdr{perso.ens-lyon.fr/leo.poyeton/}

\date{\today}

\title{Families of Galois representations and $(\varphi, \tau)$-modules}

\begin{document}

\begin{abstract}
Let $p$ be a prime, and let $K$ be a finite extension of $\Qp$, with absolute Galois group $\G_K$. Let $\pi$ be a uniformizer of $K$ and let $K_\infty$ be the Kummer extension obtained by adjoining to $K$ a system of compatible $p^n$-th roots of $\pi$, for all $n$, and let $L$ be the Galois closure of $K_\infty$. Using these extensions, Caruso has constructed étale $(\phi,\tau)$-modules, which classify $p$-adic Galois representations of $K$. In this paper, we use locally analytic vectors and theories of families of $\phi$-modules over Robba rings to prove the overconvergence of $(\phi,\tau)$-modules in families. 
As examples, we also compute some explicit families of $(\phi,\tau)$-modules in some simple cases.
\end{abstract}

\subjclass{}

\keywords{}

\maketitle

\tableofcontents

\setlength{\baselineskip}{18pt}

\section*{Introduction}

Let $K$ be a finite extension of $\Qp$. We fix an algebraic closure $\overline{K}$ of $K$ and we let $\G_K = \Gal(\overline{K}/K)$. In order to study $p$-adic Galois representations, Fontaine has constructed in \cite{Fon90} an equivalence of categories $V \mapsto D(V)$ between the category of $p$-adic representations of $\G_K$ and the category of étale $(\phi,\Gamma)$-modules. A $(\phi,\Gamma)$-module is a finite dimensional vector space over a local field $\B_K$ of dimension $2$, equipped with compatible semilinear actions of a semilinear Frobenius $\phi$ and $\Gamma$, and is said to be étale if the Frobenius is of slope $0$. In the case $K=\Qp$ (or more generally when $K$ is absolutely unramified), we can actually identify $\B_K$ to the ring of formal power series $\sum_{k \in \Z}a_kT^k$, where the sequence $(a_k)$ is a bounded sequence of elements of $K$, and such that $\lim\limits_{k \to -\infty}a_k = 0$, the actions of $\phi$ and $\Gamma_K$ being given by $\phi(T)=(1+T)^p-1$ and $\gamma(T)=(1+T)^{\chi(\gamma)}-1$, where $\chi : \G_K \to \Zp^\times$ is the cyclotomic character. 

The ring $\B_K$ does not have a satisfying analytic interpretation which makes it difficult to work with, but it contains the ring $\B_K^\dagger$ of overconvergent power series, which are the power series that converge and are bounded on an annulus bordered by the unit circle. One of the fundamental results concerning étale $(\phi,\Gamma_K)$-modules is the main theorem of \cite{cherbonnier1998representations} that shows that every étale $(\phi,\Gamma_K)$-module comes by base change from an overconvergent $(\phi,\Gamma_K)$-module defined over $\B_K^\dagger$.

Let $S$ be a $\Qp$-Banach algebra, such that for every $x$ in the maximal spectrum of $S$, $S/\mathfrak{m}_x$ is a finite extension of $\Qp$. A family of representations of $\G_K$ is a free $S$-module $V$ of finite rank, endowed with a linear continuous action of $\G_K$. In \cite{BC08}, Berger and Colmez generalized the theory of overconvergent $(\phi,\Gamma)$-modules to such families. However, in constract with the classical theory of $(\phi,\Gamma)$-modules, this functor fails to be an equivalence of categories as it is not essentially surjective. The main obstruction is the absence of a family version of Kedlaya's slope filtrations theorems, because the slope polygons of families of $\phi$-modules are not locally constant in general.

The theory of $(\phi,\Gamma)$-modules arises from the ``dévissage'' of the extension $\overline{K}/K$ through an intermediate extension $K_\infty/K$ which is the cyclotomic extension of $K$. For many reasons, one would like to generalize Fontaine's constructions by using other extensions than the cyclotomic one. The two main class of extensions considered by now are Kummer extensions (arising notably from the work of Breuil \cite{breuil1998schemas} and Kisin \cite{KisinFiso}) and Lubin-Tate extensions attached to uniformizers of a subfield of $K$, of which the cyclotomic extension is a particular case and which has been studied in order to try to extend the $p$-adic Langlands correspondence to $\GL_2(K)$ (see for exemple \cite{KR09}, \cite{FX13} or \cite{Ber14MultiLa}). In \cite{Car12}, Caruso defined the notion of $(\phi,\tau)$-modules, the analogue of $(\phi,\Gamma)$-modules for Kummer extensions. In the particular case of semi-stable representations, these $(\phi,\tau)$-modules coincide with the notion of Breuil-Kisin modules and can thus be used to study Galois deformation rings as in \cite{kisin2008potentially}, 
to classify semi-stable (integral) Galois representations as in \cite{Liu10},
and to study integral models of Shimura varieties as in \cite{kisin2010integral}. In particular, families of Breuil-Kisin modules are a particular case of families of $(\phi,\tau)$-modules.

The main goal of this paper is to construct a functor similar to the one of Berger and Colmez but for families of overconvergent étale $(\phi,\tau)$-modules. Recall that classical $(\phi,\tau)$-modules are constructed the following way: fix a Kummer extension $K_\infty=\bigcup_{n \geq 1}K_n$ where $K_n = K(\pi_n)$ and $(\pi_n)$ is a compatible sequence of $p^n$-th roots of a given uniformizer $\pi$ of $K$. As in the cyclotomic case, $p$-adic representations of $H_{\tau,K}:=\Gal(\overline{K}/K_\infty)$ are classified by étale $\phi$-modules over a local field $\B_{\tau,K}$ of dimension $2$ which can be identified with the ring of formal power series $\sum_{k \in \Z}a_kT^k$, where $(a_k)$ is a bounded sequence of elements of $F = K \cap \Q_p^{\mathrm{unr}}$, and such that $\lim\limits_{k \to -\infty}a_k = 0$. However, the comparison with the cyclotomic setting ends here, because $K_\infty/K$ is not Galois and there is thus no group action available to replace the $\Gamma$ action. The idea of Caruso is to add the action of a well chosen element $\tau$ of $\G_K$, not directly on the $\phi$-module but on the module obtained after tensoring over $\B_{\tau,K}$ by some $\B_{\tau,K}$-algebra $\Bt_L$ endowed with an action of $\Gal(L/K)$, where $L=K_\infty^{\Gal}$ (actually one can show that it's the same as adding an action of the whole group $\Gal(K_\infty^{\Gal}/K)$). As in the cyclotomic case, the ring $\B_{\tau,K}$ contains the ring $\B_{\tau,K}^\dagger$ of overconvergent elements, and one can show that that every étale $(\phi,\tau)$-module comes by base change from an overconvergent $(\phi,\tau)$-module defined over $\B_{\tau,K}^\dagger$.

The overconvergence of the classical $(\phi,\tau)$-modules has been proven by two different means in \cite{gao2016loose} and \cite{GP18} but the first proof does not extend to families. The proof used in \cite{GP18} can be generalized to construct families over the Robba ring corresponding to the Kummer extension. First, we construct a family of $\phi$-modules over $S \hat{\otimes}\Bt_{\mathrm{rig},L}$, the ``big Robba ring attached to $L$ over $S$'', endowed with an action of $\Gal(L/K)$. This can be done for example by tensoring our family of representations $V$ over the generalized Robba ring $S \hat{\otimes}\Bt_{\mathrm{rig}}$ and taking the $\Gal(\overline{K}/L)$-invariants. Then, we use the fact that we know there are ``enough'' locally analytic vectors on this level, using the results of \cite{BC08} as an input. Finally, we prove a monodromy descent theorem which allows us to descend to the level of $S\hat{\otimes}\B_{\tau,\rig,K}^\dagger$, using the computation of the pro-analytic vectors of $S\hat{\otimes}\Bt_{\tau,\rig,K}^\dagger$. The result we obtain is the following:

\begin{theo}[Theorem A]
Let $V$ be a family of representations of $\G_K$ of rank $d$. Then there exists $s_0 \geq 0$ such that for any $s \geq s_0$, there exists a unique sub-$S\hat{\otimes}\B_{\tau,\rig,K}^{\dagger,s}$-module of $(S\hat{\otimes}\Bt_{\rig,L}^{\dagger,s})^{\Gal(L/K_\infty)}$ $D_{\tau,\rig,K}^{\dagger,s}(V)$, which is a family of $(\phi,\tau)$-modules over $(S\hat{\otimes}\B_{\tau,\rig,K}^{\dagger,s},\Bt_{\rig,L}^{\dagger,s})$ such that:
\begin{enumerate}
\item $D_{\tau,\rig,K}^{\dagger,s}(V)$ is a $S\hat{\otimes}\B_{\tau,\rig,K}^{\dagger,s}$-module locally free of rank $d$;
\item the map $(S\hat{\otimes}\Bt_{\rig}^{\dagger,s})\otimes_{S\hat{\otimes}\B_{\tau,\rig,K}^{\dagger,s}}D_{\tau,\rig,K}^{\dagger,s}(V) \rightarrow (S\hat{\otimes}\Bt_{\rig}^{\dagger,s})\otimes_S V$ is an isomorphism;
\item if $x \in \cal{X}$, the map $S/\mathfrak{m}_x\otimes_SD_{\tau,\rig,K}^{\dagger,s}(V) \rightarrow D_{\tau,\rig,K}^{\dagger,s}(V_x)$ is an isomorphism.
\end{enumerate}
\end{theo}

In order to descend this family of $(\phi,\tau)$-modules over $S \hat{\otimes} \B_{\tau,\rig,K}^\dagger$ to a family of $(\phi,\tau)$-modules over the ring $S \hat{\otimes} \B_{\tau,K}^\dagger$ of bounded elements of $S \hat{\otimes} \B_{\tau,\rig,K}^\dagger$, one would like to apply a family version of Kedlaya's slope filtrations theorems. However, as explained above, there is no family version of Kedlaya's slope filtrations theorems. 

Kedlaya and Liu have proven in \cite{KL10} that the functor constructed by Berger and Colmez in \cite{BC08} could be inverted locally around rigid analytic points. The main problem is that the representations obtained this way may not glue, and the obstruction exists purely at the residual level. Considering the same question, Hellman proved in \cite{Hellmann16} that given a family $D$ of $(\phi,\Gamma)$-modules over $S\hat{\otimes}\B_{\rig,K}^\dagger$, there exists natural subfamilies $D^{\mathrm{int}}$ resp. $D^{\mathrm{adm}}$ which are étale resp. induced by a family of $p$-adic representations. Because Hellman's proof relies only on the $\phi$-structure and does not use the $\Gamma$-action, it can be translated almost directly to our setting, even though our Frobenius is not the same as the one appearing in the cyclotomic theory. Because we already know the family of $\phi$-modules $D$ over the Robba ring obtained in theorem A comes from a family of $p$-adic representations, this implies that we have $D = D^{\mathrm{int}} = D^{\mathrm{adm}}$ and allows us to recover the family of $(\phi,\tau)$-modules over the ring $S \hat{\otimes} \B_{\tau,K}^\dagger$.

\begin{theo}[Theorem B]
Let $V$ be a family of representations of $\G_K$ of rank $d$. Then there exists $s_0 \geq 0$ such that for any $s \geq s_0$, there exists a family of $(\phi,\tau)$-modules $D_{\tau,K}^{\dagger,s}(V)$ such that:
\begin{enumerate}
\item $D_{\tau,K}^{\dagger,s}(V)$ is a $S\hat{\otimes}\B_{\tau,K}^{\dagger,s}$-module locally free of rank $d$;
\item the map $(S\hat{\otimes}\Bt^{\dagger,s})\otimes_{S\hat{\otimes}\B_{\tau,K}^{\dagger,s}}D_{\tau,K}^{\dagger,s}(V) \rightarrow (S\hat{\otimes}\Bt^{\dagger,s})\otimes_S V$ is an isomorphism;
\item if $x \in \cal{X}$, the map $S/\mathfrak{m}_x\otimes_SD_{\tau,K}^{\dagger,s}(V) \rightarrow D_{\tau,K}^{\dagger,s}(V_x)$ is an isomorphism.
\end{enumerate}
\end{theo}

As an example, we then compute the family of all rank $1$ $(\phi,\tau)$-modules in the particular case $K=\Qp$ and $\pi=p$. Let $T$ be such that $\B_{\tau,\Qp}$ is the ring of formal power series $\sum_{k\in \Z}a_kT^k$, with $\phi(T)=T^p$, let $\lambda = \prod_{k=0}^\infty \phi^n(\frac{T-p}{p})$. We let $\mu_{\beta}$ denote the character of $\Q_p^{\times}$ sending $p$ to $\beta$ and which is trivial on $\Z_p^\times$. 

\begin{theo}[Theorem C]
There exists $\alpha \in \At_L^+$ such that the $(\phi,\tau)$-module corresponding to $\delta = \mu_{\beta}\cdot \omega^r \cdot \langle \chi_{\cycl} \rangle^s$ admits a basis $e$ in which $\phi(e) = \beta \cdot T^r\cdot (1-\frac{p}{T})^{-s}\cdot e$ and $\tau(e) = [\epsilon]^{-r}\prod_{n=0}^{+\infty}\phi^n(1+\alpha T)^{-s}\cdot e$.
\end{theo}

As an other example, we give a description of the Breuil-Kisin modules attached to some trianguline semistable representations of rank $2$. The notion of trianguline representations was introduced by Colmez in \cite{colmez2010representations}, which are representations whose attached $(\phi,\Gamma)$-module over the Robba ring is a successive extension of rank $1$ $(\phi,\Gamma)$-modules. First, we show that this is equivalent for a representation to be trianguline in the sense of Colmez, and to be trianguline in the sense of $(\phi,\tau)$-modules, that is that its $(\phi,\tau)$-module is a successive extension of rank $1$ $(\phi,\tau)$-modules. Using this and the fact that Breuil-Kisin modules can be constructed directly by using $(\phi,\tau)$-modules, we show the following, which is a consequence of theorem C, of the compatibility between $(\phi,\tau)$-modules and Breuil-Kisin modules, and of Kisin's results \cite{KisinFiso}. Note that in what follows we let $\lambda '=\frac{d}{dT}\lambda$. (Please see Section \ref{sec:Explicit} for precise notation.) 

\begin{theo}[Theorem D]
Let $V$ be a trianguline semistable representation, with nonpositive Hodge-Tate weights, whose $(\phi,\Gamma)$-module is an extension of $\cal{R}(\delta_1)$ by $\cal{R}(\delta_2)$, where $\delta_1(\Z_p^\times)$ and $\delta_2(\Z_p^\times)$ belong to $\Z_p^\times$, and are respectively of weight $k_1$ and $k_2$. Then the $(\phi,N_\nabla)$-module attached to $V$ admits a basis in which
\begin{equation*}
\Mat(\phi)= 
\begin{pmatrix}
\delta_1(p)(T-p)^{-k_1} & (T-p)^{\inf(-k_1,-k_2)}\alpha_V \\
0 & \delta_2(p)(T-p)^{-k_2} 
\end{pmatrix}
\end{equation*}
and
\begin{equation*}
\Mat(N_\nabla)= 
\begin{pmatrix}
-k_1T\lambda' & \beta_V \\
0 & -k_2T\lambda' 
\end{pmatrix},
\end{equation*}
where $\alpha_V, \beta_V \in \B_{\tau,\rig,K}^\dagger$. Moreover, $V$ is crystalline if and only if $\beta_V = 0 \mod [\tilde{p}]$. 
\end{theo}

\subsection*{Organization of the paper}
This paper is subdivided into 7 sections. The first one is devoted to the reminders on Kummer extensions and the definitions of $(\phi,\tau)$-modules and $(\phi,\Gamma)$-modules and the relevant period rings that appear in the theory. The second sections recalls the basics of the theory of locally analytic vectors, its specialization in our case and the properties we shall need for the rest of the paper. In section 3, we recall a gluing property for coherent sheaves on noetherian adic spaces. In section 4, we define families of representations and recall the main results from \cite{BC08} that attach to such families corresponding families of $(\phi,\Gamma)$-modules. Section 5 shows how to construct a family of $(\phi,\tau)$-modules over the corresponding Robba ring, starting with a family of $(\phi,\Gamma)$-modules over the ``cyclotomic Robba ring'', which can then be applied to give theorem A. In section 6, we show how to descend such a family attached to a family of representations, and constructed by using Berger's and Colmez's results as an input, to a family of $(\phi,\tau)$-modules over the bounded Robba ring, which is our theorem B. Finally, the explicit computations of $(\phi,\tau)$-modules are done in the last section, in which theorems C and D appear.

\subsection*{Acknowledgements} 
The authors would like to thank Ruochuan Liu for helpful conversations and comments.

\section{Kummer extensions and $(\phi,\tau)$-modules}
\subsection{Kummer extensions and first definitions}
Let $K$ be a finite extension of $\Qp$, with ring of integers $\O_K$ and residue field $k$ and let $\pi$ be a uniformizer of $K$. Let $F=W(k)[1/p]$, so that $K/F$ is a finite totally ramified extension, and let $e=[K:F]$. Let $E(X) \in \O_F[X]$ be the minimal poynomial of $\pi$ over $F$. We let $\overline{K}$ be an algebraic closure of $K$ and let $\C_p$ be the $p$-adic completion of $\overline{K}$. Let $v_p$ be the $p$-adic valuation on $\Cp$ such that $v_p(p)=1$. Let $(\pi_n)$ be a sequence of elements of $\overline{K}$, such that $\pi_0 =\pi$ and $\pi_{n+1}^p = \pi_n$. We let $K_n = K(\pi_n)$ and $K_\infty = \bigcup_{n \geq 1}K_n$. Let also $\epsilon_1$ be a primitive $p$-th root of unity and $(\epsilon_n)_{n \in \N}$ be a compatible sequence of $p^n$-th roots of unity, which means that $\epsilon_{n+1}^p=\epsilon_n$ and let $K_{\mathrm{cycl}} = \bigcup_{n \geq 0}K(\epsilon_n)$ be the cyclotomic extension of $K$. Let $L:=K_\infty \cdot K_{\mathrm{cycl}}$ be the Galois closure of $K_\infty/K$, and let 
\[
G_\infty = \Gal(L/K), \quad H_\infty = \G_L = \Gal(\overline{K}/L), \quad \Gamma_K = \Gal(L/K_\infty).
\]
Note that we can identify $\Gamma_K$ with $\Gal(K_{\mathrm{cycl}}/(K_\infty \cap K_{\mathrm{cycl}}))$ and so to an open subgroup of $\Z_p^\times$. 

For $g \in \G_K$ and for $n \in \N$, there exists a unique element $c_n(g) \in \Z/p^n\Z$ such that $g(\pi_n) = \epsilon_n^{c_n(g)}\pi_n$. Since $c_{n+1}(g) = c_n(g) \mod p^n$, the sequence $(c_n(g))$ defines an element $c(g)$ of $\Zp$. The map $g \mapsto c(g)$ is actually a (continuous) $1$-cocycle of $\G_K$ to $\Zp(1)$, such that $c^{-1}(0) = \Gal(\overline{K}/K_{\infty})$, and satisfies for $g,h \in \Gal(\overline{K}/K_{\infty})$~:
$$c(gh) = c(g)+\chi_{\mathrm{cycl}}(g)c(h).$$

This means that if $\Zp \rtimes \Z_p^{\times}$ is the semi-direct product of $\Zp$ by $\Z_p^{\times}$ where $\Z_p^{\times}$ acts on $\Zp$ by multiplication, then the map $g \in \G_K \mapsto (c(g),\chi_{\mathrm{cycl}}(g)) \in \Zp \rtimes \Z_p^{\times}$ is a morphism of groups of Kernel $H_{\infty}$. The cocycle $c$ factors through $H_{\infty}$, and so defines a cocycle that we will still denote by $c~: G_\infty \to \Zp$ which is the Kummer's cocycle attached to $K_\infty/K$.

We let $\tau$ be a topological generator of $\Gal(L/K_{\mathrm{cycl}})$ such that $c(\tau)=1$ (this is exactly the element corresponding to $(1,1)$ \textit{via} the isomorphism $g \in \G_L \mapsto (c(g),\chi_{\mathrm{cycl}}(g)) \in \Zp \rtimes \Z_p^\times$). The relation between $\tau$ and $\Gamma_K$ is given by $g\tau g^{-1} = \tau^{\chi_{\mathrm{cycl}}(g)}$. We also let $\tau_n:=\tau^{p^n}$.

Since we will consider at the same time rings relative to the cyclotomic extension of $K$ and rings relative to the Kummer extension $K_\infty$ of $K$, we will write a $\tau$ in index of the rings relative to the Kummer extension. Note that it does not depend on the choice of $\tau$. We also let $H_K = \Gal(\overline{K}/K_{\mathrm{cycl}})$ and $H_{\tau,K} = \Gal(\overline{K}/K_\infty)$. If $A$ is an algebra endowed with an action of $\G_K$, we let $A_K = A^{H_K}$ and $A_{\tau,K} = A^{H_{\tau,K}}$.

\subsection{$(\phi,\tau)$-modules, $(\phi,\Gamma)$-modules and some (most?) involved rings}
Let
\[\Etplus = \varprojlim\limits_{x \to x^p} \O_{\C_p} = \{(x^{(0)},x^{(1)},\dots) \in \O_{\C_p}^{\N}~: (x^{(n+1)})^p=x^{(n)}\}
\]
and recall that $\Etplus$ is naturally endowed with a ring structure that makes it a perfect ring of characteristic $p$ which is complete for the valuation $v_{\E}$ defined by $v_{\E}(x) = v_p(x^{(0)})$. Let $\Et$ be its fraction field. The theory of field of norms of Fontaine-Wintenberger \cite{Win83} allows us to attach to the extension $K_\infty/K$ its field of norms $X_K(K_\infty)$ which injects into $\Et$. The sequences $(\epsilon_n)$ and $(\pi_n)$ define elements of $\Etplus$ which we will denote respectively by $\epsilon$ and $\tilde{\pi}$. Let $\overline{u} = \epsilon -1$, and recall that $v_{\E}(\overline{u}) = \frac{p}{p-1}$. The image of the injection of $X_K(K_\infty)$ inside $\Et$ is then $\E_{\tau,K} := k(\!(\tilde{\pi})\!)$. Let $\E_\tau$ be the separable closure of $\E_{\tau,K}$ inside $\Et$. Since $\Gal(\overline{K}/K_{\infty})$ acts trivially on $\E_{\tau,K}$, every element of $\Gal(\overline{K}/K_{\infty})$ stabilizes $\E_\tau$, which gives us a morphism $\Gal(\overline{K}/K_{\infty}) \to \Gal(\E_\tau/\E_{\tau,K})$ which is an isomorphism by theorem 3.2.2 of \cite{Win83}.

Let 
\[
\At = W(\Et),  \quad\Atplus = W(\Etplus), \quad \Bt = \At[1/p] \quad \textrm{and } \Btplus = \Atplus[1/p].
\]
These rings are equipped with a Frobenius $\phi$ deduced from the one on $\Et$ by the functoriality of Witt vectors and with a $\G_\Qp$-action lifting the one on $\Et$ and given by $g\cdot [x]=[g\cdot x]$.

These rings are naturally endowed with two different topologies, called respectively the strong topology and the weak topology. The strong topology is the coarsest topology such that the projection $\At \to \Et$ is continuous, where $\Et$ is endowed with the discrete topology. Hence, the strong topology on $\At$ is the $p$-adic topology. The weak topology is the coarsest topology such that the projection $\At \to \Et$ is continuous, where $\Et$ is endowed with the topology given by $v_{\E}$. By \cite[Prop. 5.2]{colmez2008espaces}, the action of $\phi$ on $\At$, $\Atplus$, $\Bt$ and $\Btplus$ is continuous for both the strong and the weak topology, and we have 
$$\At^{\phi=1} = (\Atplus)^{\phi=1} = \Zp, \quad \Bt^{\phi=1} = (\Btplus)^{\phi=1} = \Qp.$$

Since neither $H_K$, $H_{\tau,K}$ nor $\G_\Qp$ act continuously on $\Et$ or $\Etplus$ for the discrete topology, they don't act continuously on $\At$ or $\Atplus$ for the $p$-adic topology. However, since $\At$ and $\At^+$ are respectively homeomorphic to $\Et^{\N}$ and $(\Et^+)^{\N}$, $\G_\Qp$ acts continuously on $\At$ and $\At^+$ endowed with the weak topology.

For $r > 0$, we define $\Bt^{\dagger,r}$ the subset of overconvergent elements of ``radius'' $r$ of $\Bt$, by

$$\left\{x = \sum_{n \ll -\infty}p^n[x_n] \textrm{ such that } \lim\limits_{k \to +\infty}v_{\E}(x_k)+\frac{pr}{p-1}k =+\infty \right\}$$

and we let $\Bt^\dagger = \bigcup_{r > 0}\Bt^{\dagger,r}$ be the subset of all overconvergent of $\Bt$. 

We now define a ring $\A_{\tau,K}$ inside $\At$ as follows:
$$\A_{\tau,K} = \left\{\sum_{i \in \Z}a_i[\tilde{\pi}]^i, a_i \in \O_F \lim\limits_{i \to - \infty}a_i = 0 \right\}.$$
Endowed with the $p$-adic valuation $v_p(\sum_{i \in \Z}a_i[\tilde{\pi}]^i) = \min_{i \in \Z}v_p(a_i)$, $\A_{\tau,K}$ is a DVR with residue field $\E_{\tau,K}$. Let $\B_{\tau,K}=\A_{\tau,K}[1/p]$ and let $\B_{\tau,K}^{\dagger,r}$ be the subset of $\B_{\tau,K}$ given by
$$\B_{\tau,K}^{\dagger,r}=\left\{\sum_{i \in \Z}a_i[\tilde{\pi}]^i, a_i \in F \textrm{ such that the } a_i \textrm{ are bounded and } \lim\limits_{i \to - \infty}v_p(a_i)+i\frac{pr}{p-1} = +\infty \right\}.$$
Note that $\B_{\tau,K}^{\dagger,r} = \B_{\tau,K} \cap \Bt^{\dagger,r}$.

Let $\B_{\tau,K}^\dagger = \bigcup_{r > 0}\B_{\tau,K}^{\dagger,r}$. By §2 of \cite{matsuda1995local}, this is a Henselian field, and its residue ring is still $\E_{\tau,K}$. If $M/K$ is a finite extension, we let $\E_{\tau,M}$ be the extension of $\E_{\tau,K}$ corresponding to $M\cdot K_\infty/K_\infty$ by the theory of field of norms, which is a separable extension of degree $f=[M\cdot K_\infty:K_\infty]$. Since $\B_{\tau,K}^\dagger$ is Henselian, there exists a finite unramified extension $\B_{\tau,M}^\dagger/\B_{\tau,K}^\dagger$ inside $\Bt$, of degree $f$ and whose residue field is $\E_{\tau,M}$. Therefore, there exists $r(M) > 0$ and elements $x_1,\ldots,x_f$ in $\B_{\tau,M}^{\dagger,r(M)}$ such that $\B_{\tau,M}^{\dagger,s} = \bigoplus_{i=1}^f \B_{\tau,K}^{\dagger,s}\cdot x_i$ for all $s \geq r(M)$. We let $\B_{\tau,M}$ be the $p$-adic completion of $\B_{\tau,M}^\dagger$. 

The Frobenius on $\Bt$ defines by restriction endomorphisms of $\A_{\tau,K}$ and $\B_{\tau,K}$, and sends $[\tilde{\pi}]$ to $[\tilde{\pi}]^p$. We also let $\At_L = \At^{H_{\infty}}$ and $\Bt_L = \At_L[1/p]$.

A $\phi$-module $D$ on $\B_{\tau,K}$ is a $\B_{\tau,K}$-vector space of finite dimension $d$, equipped with a semilinear $\phi$ action such that $\Mat(\phi) \in \GL_d(\B_{\tau,K})$, and we say that it is étale if there exists a basis of $D$ in which $\Mat(\phi) \in \GL_d(\A_{\tau,K})$.

Usual $(\phi,\tau)$-modules can be defined as follows:
\begin{defi}
\label{def phitau}
A $(\phi,\tau)$-module on $(\B_{\tau,K},\Bt_L)$ is a triple $(D,\phi_D,G)$ where:
\begin{enumerate}
\item $(D,\phi_D)$ is a $\phi$-module on $\B_{\tau,K}$;
\item $G$ is a continuous (for the weak topology) $G_\infty$-semilinear $G_\infty$-action on $M:=\Bt_L \otimes_{\B_{\tau,K}}D$ such that $G$ commutes with $\phi_M:=\phi_{\Bt_L}\otimes \phi_D$, i.e. for all $g \in G_\infty$, $g\phi_M = \phi_Mg$;
\item regarding $D$ as a sub-$\B_{\tau,K}$-module of $M$, $D \subset M^{H_{\tau,K}}$.
\end{enumerate}
We say that a $(\phi,\tau)$-module is étale if its underlying $\phi$-module on $\B_{\tau,K}$ is.
\end{defi}
This definition is the same as \cite[Def. 2.1.5]{gao2016loose} and not the same as Caruso's, however note that both definitions are equivalent for $p \neq 2$ (see remark 2.1.6 of \cite{gao2016loose}) and that this definition ``works'' in the case $p=2$, meaning that we have the following:

\begin{prop}
\label{prop eqcat etalephitau padicrep}
Given an étale $(\phi,\tau)$-module $(D,\phi_D,G)$, we define 
$$V(D):=(\Bt \otimes_{\Bt_L}M)^{\phi=1},$$
where $M=\Bt_L \otimes_{\B_{\tau,K}}D$ is equipped with a $\G_\infty$-action. Note that $V(D)$ is a $\Qp$-vector space endowed with a $\G_K$ action induced by the ones on $\Bt$ and $M$.

The functor $D \mapsto V(D)$ induces an equivalence of categories between the category of étale $(\phi,\tau)$-modules and the category of $p$-adic representations of $\G_K$.
\end{prop}
\begin{proof}
This is \cite[Prop. 2.1.7]{gao2016loose}.
\end{proof}

We also quickly recall some of the theory of $(\phi,\Gamma)$-modules and the definitions of some rings that appear in this theory, as we will need them later on.

Let $\E_{\Qp} = \F_p(\!(\overline{u})\!) \subset \Et$. This is the image of the field of norms $X_{\Qp}((\Q_p)_{\mathrm{cycl}})$. Let $u = [\epsilon]-1 \in \Atplus$. We define a ring $\A_{\Qp}$ inside $\At$ as follows:
$$\A_{\Qp} = \left\{\sum_{i \in \Z}a_iu^i, a_i \in \Zp \lim\limits_{i \to - \infty}a_i = 0 \right\}.$$
Endowed with the $p$-adic valuation, $\A_{\Qp}$ is a DVR with residue field $\E_{\Qp}$. Let $\B_{\Qp}=\A_{\Qp}[1/p]$. The Frobenius on $\Bt$ defines by restriction an endomorphism on $\B_{\Qp}$, and we also have an action of $\G_{\Qp}$ on $\B_{\Qp}$. These actions are given by
$$\phi(u)=(1+u)^p-1, \quad g(u) = (1+u)^{\chi_{\cycl}(g)}-1$$
and commute with each other.

Let $\B_{\Qp}^{\dagger,r}$ be the subset of $\B_{\Qp}$ given by
$$\B_{\Qp}^{\dagger,r}=\left\{\sum_{i \in \Z}a_i[\tilde{\pi}]^i, a_i \in \Qp \textrm{ such that the } a_i \textrm{ are bounded and } \lim\limits_{i \to - \infty}v_p(a_i)+i\frac{pr}{p-1} = +\infty \right\},$$
and note that $\B_{\Qp}^{\dagger,r} = \B_{\Qp} \cap \Bt^{\dagger,r}$.  

Let $\B_{\Qp}^\dagger = \bigcup_{r > 0}\B_{\Qp}^{\dagger,r}$. By §2 of \cite{matsuda1995local}, this is a Henselian field, and its residue ring is still $\E_{\Qp}$. If $M/\Qp$ is a finite extension, we let $\E_M$ be the extension of $\E_{\Qp}$ corresponding to $M_{\cycl}/\Qp_{\cycl}$ by the theory of field of norms, which is a separable extension of degree $f=[M_{\cycl}:\Qp_{\cycl}]$. Since $\B_{\Qp}^\dagger$ is Henselian, there exists a finite unramified extension $\B_M^\dagger/\B_{\Qp}^\dagger$ inside $\Bt$, of degree $f$ and whose residue field is $\E_M$. Therefore, there exists $r(M) > 0$ and elements $x_1,\ldots,x_f$ in $\B_M^{\dagger,r(M)}$ such that $\B_M^{\dagger,s} = \bigoplus_{i=1}^f \B_{\Qp}^{\dagger,s}\cdot x_i$ for all $s \geq r(M)$. We let $\B_M$ be the $p$-adic completion of $\B_M^\dagger$ and we let $\A_M$ be its ring of integers for the $p$-adic valuation. One can show that $\B_M$ is a subfield of $\Bt$ stable under the action of $\phi$ and $\Gamma_M = \Gal(M_{\cycl}/M)$ (see for example \cite[Prop. 6.1]{colmez2008espaces}).

A $\phi$-module $D$ on $\B_K$ is a $\B_{K}$-vector space of finite dimension $d$, equipped with a semilinear $\phi$ action such that $\Mat(\phi) \in \GL_d(\B_K)$, and we say that it is étale if there exists a basis of $D$ in which $\Mat(\phi) \in \GL_d(\A_K)$.
We can now define the notion of $(\phi,\Gamma_K)$-modules on $\B_K$:
\begin{defi}
A $(\phi,\Gamma_K)$-module $D$ on $\B_K$ is a $\phi$-module on $\B_K$ equipped with a commuting and continuous (for the weak topology) semilinear action of $\Gamma_K$. We say that it is étale if the underlying $\phi$-module is.
\end{defi}

We then have the following proposition:
\begin{prop}
Given an étale $(\phi,\Gamma_K)$-module $D$, we define
$$V(D):=(\Bt \otimes_{\B_{\tau,K}}D)^{\phi=1},$$
which is a $\Qp$-vector space endowed with a $\G_K$ action coming from the ones on $\Bt$ and $D$.

The functor $D \mapsto V(D)$ induces an equivalence of categories between the category of étale $(\phi,\Gamma_K)$-modules on $\B_{\tau,K}$ and the category of $p$-adic representations of $\G_K$.
\end{prop}
\begin{proof}
This is \cite[Thm. A.3.4.3]{Fon90}.
\end{proof}

\subsection{Some more rings of periods}
For $r \geq 0$, we define a valuation $V(\cdot,r)$ on $\Btplus[\frac{1}{[\tilde{\pi}]}]$ by setting
$$V(x,r) = \inf_{k \in \Z}(k+\frac{p-1}{pr}v_{\E}(x_k))$$
for $x = \sum_{k \gg - \infty}p^k[x_k]$. If $I$ is a closed subinterval of $[0;+\infty[$, we let $V(x,I) = \inf_{r \in I}V(x,r)$. We then define the ring $\Bt^I$ as the completion of $\Btplus[1/[\tilde{\pi}]]$ for the valuation $V(\cdot,I)$ if $0 \not \in I$, and as the completion of $\Btplus$ for $V(\cdot,I)$ if $I=[0;r]$. We will write $\Bt_{\mathrm{rig}}^{\dagger,r}$ for $\Bt^{[r,+\infty[}$ and $\Bt_{\mathrm{rig}}^+$ for $\Bt^{[0,+\infty[}$. We also define $\Bt_{\mathrm{rig}}^\dagger = \bigcup_{r \geq 0}\Bt_{\mathrm{rig}}^{\dagger,r}$.

Let $I$ be a subinterval of $]1,+\infty[$ or such that $0 \in I$. Let $f(Y) = \sum_{k \in \Z}a_kY^k$ be a power series with $a_k \in F$ and such that $v_p(a_k)+\frac{p-1}{pe}k/\rho \to +\infty$ when $|k| \to + \infty$ for all $\rho \in I$. The series $f([\tilde{\pi}])$ converges in $\Bt^I$ and we let $\B_{\tau,K}^I$ denote the set of all $f([\tilde{\pi}])$ with $f$ as above. It is a subring of $\Bt_{\tau,K}^I$. The Frobenius gives rise to a map $\phi: \B_{\tau,K}^I \to \B_{\tau,K}^{pI}$. If $m \geq 0$, then we have $\phi^{-m}(\B_{\tau,K}^{p^mI}) \subset \Bt_{\tau,K}^I$ and we let $\B_{\tau,K,m}^I = \phi^{-m}(\B_{\tau,K}^{p^mI})$, so that $\B_{\tau,K,m}^I \subset \B_{\tau,K,m+1}^I$ for all $m \geq 0$.

We also write $\B_{\tau,\mathrm{rig},K}^{\dagger,r}$ for $\B_{\tau,K}^{[r;+\infty[}$. It is a subring of $\B_{\tau,K}^{[r;s]}$ for all $s \geq r$ and note that the set of all $f([\tilde{\pi}]) \in \B_{\tau,\mathrm{rig},K}^{\dagger,r}$ such that the sequence $(a_k)_{k \in \Z}$ is bounded is exactly the ring $\B_{\tau,K}^{\dagger,r}$. Let $\B_{\tau,K}^{\dagger}=\cup_{r \gg 0}\B_{\tau,K}^{\dagger,r}$. Let $\B_{\tau,K,m}^I = \phi^{-m}(\B_{\tau,K}^{p^mI})$ and $\B_{\tau,K,\infty}^I=\cup_{m \geq 0}\B_{\tau,K,m}^I$ so that in particular we have $\B_{\tau,K,m}^I \subset \Bt_{\tau,K}^I$. 

Recall that, for $M$ a finite extension of $K$, there exists by the theory of field of norms a separable extension $\E_{\tau,M}/\E_{\tau,K}$ of degree $f=[M_{\infty}:K_{\infty}]$ and an attached unramified extension $\B_{\tau,M}^{\dagger}/\B_{\tau,K}^{\dagger}$ of degree $f$ with residue field $\E_{\tau,M}$, so that there exists $r(M) > 0$ and elements $x_1,\cdots x_f \in \B_{\tau,M}^{\dagger,r(M)}$ such that $\B_{\tau,M}^{\dagger,s}= \bigoplus_{i=1}^f\B_{\tau,K}^{\dagger,s}\cdot x_i$ for all $s \geq r(M)$. If $r(M) \leq \min(I)$, we let $\B_{\tau,M}^I$ be the completion of $\B_{\tau,M}^{\dagger,r(M)}$ for $V(\cdot,I)$, so that $\B_{\tau,M}^I=\oplus_{i=1}^f\B_{\tau,K}^I\cdot x_i$. 

We will also define the corresponding rings for the cyclotomic setting.

Let $I$ be a subinterval of $]1,+\infty[$ or such that $0 \in I$. Let $f(Y) = \sum_{k \in \Z}a_kY^k$ be a power series with $a_k \in F$ and such that $v_p(a_k)+\k/\rho \to +\infty$ when $|k| \to + \infty$ for all $\rho \in I$. The series $f(u)$ converges in $\Bt^I$ and we let $\B_{\Qp}^I$ denote the set of all $f(u)$ with $f$ as above. It is a subring of $\Bt_{\Qp}^I$. The Frobenius gives rise to a map $\phi: \B_{\Qp}^I \to \B_{\Qp}^{pI}$. If $m \geq 0$, then we have $\phi^{-m}(\B_{\Qp}^{p^mI}) \subset \Bt_{\Qp}^I$ and we let $\B_{\Qp,m}^I = \phi^{-m}(\B_{\Qp}^{p^mI})$, so that $\B_{\Qp,m}^I \subset \B_{\Qp,m+1}^I$ for all $m \geq 0$.

We also write $\B_{\mathrm{rig},\Qp}^{\dagger,r}$ for $\B_{\Qp}^{[r;+\infty[}$. It is a subring of $\B_{\Qp}^{[r;s]}$ for all $s \geq r$ and note that the set of all $f(u) \in \B_{\mathrm{rig},\Qp}^{\dagger,r}$ such that the sequence $(a_k)_{k \in \Z}$ is bounded is exactly the ring $\B_{\Qp}^{\dagger,r}$. Let $\B_{\Qp}^{\dagger}=\cup_{r \gg 0}\B_{\Qp}^{\dagger,r}$. Let $\B_{\Qp,m}^I = \phi^{-m}(\B_{\Qp}^{p^mI})$ and $\B_{\Qp,\infty}^I=\cup_{m \geq 0}\B_{\Qp,m}^I$ so that in particular we have $\B_{\Qp,m}^I \subset \Bt_{\Qp}^I$. 

Recall that, for $M$ a finite extension of $\Qp$, there exists by the theory of field of norms a separable extension $\E_{M}/\E_{\Qp}$ of degree $f=[M_{\cycl}:(\Qp)_{\cycl}]$ and an attached unramified extension $\B_{M}^{\dagger}/\B_{\Qp}^{\dagger}$ of degree $f$ with residue field $\E_{M}$, so that there exists $r(M) > 0$ and elements $x_1,\cdots x_f \in \B_{M}^{\dagger,r(M)}$ such that $\B_{M}^{\dagger,s}= \bigoplus_{i=1}^f\B_{\Qp}^{\dagger,s}\cdot x_i$ for all $s \geq r(M)$. If $r(M) \leq \min(I)$, we let $\B_{M}^I$ be the completion of $\B_{M}^{\dagger,r(M)}$ for $V(\cdot,I)$, so that $\B_{M}^I=\oplus_{i=1}^f\B_{\Qp}^I\cdot x_i$.

\section{Locally analytic and pro-analytic vectors}
\subsection{Basics of the theory and key lemmas}
In this section, we recall the theory of locally analytic vectors of Schneider and Teitelbaum \cite{schneider2002bis} but here we follow the constructions of Emerton \cite{emerton2004locally} as in \cite{Ber14MultiLa}. In this article, we will use the following multi-index notations: if $\cbf = (c_1, \hdots,c_d)$ and $\kbf = (k_1,\hdots,k_d) \in \N^d$ (here $\N=\Z^{\geq 0}$), then we let $\cbf^\kbf = c_1^{k_1} \cdot \ldots \cdot c_d^{k_d}$.

Let $G$ be a $p$-adic Lie group, and let $W$ be a $\Qp$-Banach representation of $G$. Let $H$ be an open subgroup of $G$ such that there exists coordinates $c_1,\cdots,c_d : H \to \Zp$ giving rise to an analytic bijection $\cbf : H \to \Z_p^d$. We say that $w \in W$ is an $H$-analytic vector if there exists a sequence $\left\{w_{\kbf}\right\}_{\kbf \in \N^d}$ such that $w_{\kbf} \rightarrow 0$ in $W$ and such that $g(w) = \sum_{\kbf \in \N^d}\cbf(g)^{\kbf}w_{\kbf}$ for all $g \in H$. We let $W^{H-\an}$ be the space of $H$-analytic vectors. This space injects into $\cal{C}^{\an}(H,W)$, the space of all analytic functions $f : H \to W$.  Note that $\cal{C}^{\an}(H,W)$ is a Banach space equipped with its usual Banach norm, so that we can endow $W^{H-\an}$ with the induced norm, that we will denote by $||\cdot ||_H$. With this definition, we have $||w||_H = \sup_{\kbf \in \N^d}||w_{\kbf}||$ and $(W^{H-\an},||\cdot||_H)$ is a Banach space.

The space $\cal{C}^{\an}(H,W)$ is endowed by an action of $H \times H \times H$, given by
\[
((g_1,g_2,g_3)\cdot f)(g) = g_1 \cdot f(g_2^{-1}gg_3)
\]
and one can recover $W^{H-\an}$ as the closed subspace of $\cal{C}^{\an}(H,W)$ of its $\Delta_{1,2}(H)$-invariants,  where $\Delta_{1,2} : H \to H \times H \times H$ denotes the map $g \mapsto (g,g,1)$ (we refer the reader to \cite[§3.3]{emerton2004locally} for more details).

We say that a vector $w$ of $W$ is locally analytic if there exists an open subgroup $H$ as above such that $w \in W^{H-\an}$. Let $W^{\la}$ be the space of such vectors, so that $W^{\la} = \bigcup_{H}W^{H-\an}$, where $H$ runs through a sequence of open subgroups of $G$. The space $W^{\la}$ is naturally endowed with the inductive limit topology, so that it is an LB space. 

It is often useful to work with a set of ``compatible coordinates'' of $G$, that is taking an open compact subgroup $G_1$ of $G$ such that there exists local coordinates $\cbf : G_1 \to (\Zp)^d$ such that, if $G_n = G_1^{p^{n-1}}$ for $n \geq 1$, then $G_n$ is a subgroup of $G_1$ satisfying $\cbf(G_n) = (p^n\Zp)^d$. By the discussion following \cite[Prop. 2.3]{Ber14SenLa}, it is always possible to find such a subgroup $G_1$ (note however that it is not unique). We then have $W^{\la} = \bigcup_{n \in \N}W^{G_n-\an}$.

In the rest of this article, we will need the following results, most of which appear in \cite[§2.1]{Ber14SenLa} or \cite[§2]{Ber14MultiLa}.

\begin{lemm}
\label{Gn-an subset Gm-an}
Let $G_1$ and $(G_n)_{n \in \N}$ be as in the discussion above. Suppose $w \in W^{G_n-\an}$. Then for all $m \geq n$, we have $w \in W^{G_m-\an}$ and $||w||_{G_m} \leq ||w||_{G_n}$. Moreover, we have $||w||_{G_m} = ||w||$ when $m \gg 0$. 
\end{lemm}
\begin{proof}
This is \cite[Lemm. 2.4]{Ber14SenLa}.
\end{proof}

\begin{lemm}
\label{ringla}
If $W$ is a ring  such that $||xy|| \leq ||x|| \cdot ||y||$ for $x,y \in W$, then
\begin{enumerate}
  \item $W^{H-\an}$ is a ring, and $||xy||_H \leq||x||_H \cdot ||y||_H$ if $x,y \in W^{H-\an}$;
  \item if $w \in W^\times cap W^{\la}$, then $1/w \in W^{\la}$. In particular, if $W$ is a field, then  $W^{\la}$ is also a field.
\end{enumerate}
\end{lemm}
\begin{proof}
See \cite[Lemm. 2.5]{Ber14SenLa}.
\end{proof}

Let $W$ be a Fréchet space whose topology is defined by a sequence $\left\{p_i\right\}_{i \geq 1}$ of seminorms. Let $W_i$ be the Hausdorff completion of $W$ at $p_i$, so that $W = \varprojlim\limits_{i \geq 1}W_i$. The space $W^{\la}$ can be defined but as stated in \cite{Ber14MultiLa}, this space is too small in general for what we are interested in, and so we make the following definition, following \cite[Def. 2.3]{Ber14MultiLa}:

\begin{defi}
If $W = \varprojlim\limits_{i \geq 1}W_i$ is a Fréchet representation of $G$, then we say that a vector $w \in W$ is pro-analytic if its image $\pi_i(w)$ in $W_i$ is locally analytic for all $i$. We let $W^{\pa}$ denote the set of all pro-analytic vectors of $W$. 
\end{defi}

We extend the definition of $W^{\la}$ and $W^{\pa}$ for LB and LF spaces respectively. 

\begin{prop}
\label{lainla and painpa}
Let $G$ be a $p$-adic Lie group, let $B$ be a Banach $G$-ring and let $W$ be a free $B$-module of finite rank, equipped with a compatible $G$-action. If the $B$-module $W$ has a basis $w_1,\ldots,w_d$ in which $g \mapsto \Mat(g)$ is a globally analytic function $G \to \GL_d(B) \subset M_d(B)$, then
\begin{enumerate}
\item $W^{H-\an} = \bigoplus_{j=1}^dB^{H-\an}\cdot w_j$ if $H$ is a subgroup of $G$;
\item $W^{\la} = \bigoplus_{j=1}^dB^{\la}\cdot w_j$.
\end{enumerate}
Let $G$ be a $p$-adic Lie group, let $B$ be a Fréchet $G$-ring and let $W$ be a free $B$-module of finite rank, equipped with a compatible $G$-action. If the $B$-module $W$ has a basis $w_1,\ldots,w_d$ in which $g \mapsto \Mat(g)$ is a pro-analytic function $G \to \GL_d(B) \subset M_d(B)$, then
$$W^{\pa} = \bigoplus_{j=1}^dB^{\pa}\cdot w_j.$$
\end{prop}
\begin{proof}
The part for Banach ring is proven in \cite[Prop. 2.3]{Ber14SenLa} and the one for Fréchet rings is proven in \cite[Prop. 2.4]{Ber14MultiLa}.
\end{proof}

\begin{prop}
\label{prop trivial action = standard loc ana}
Let $V$ and $W$ be two $\Qp$-Banach representations of $G$ and assume that $G$ acts trivially on $W$. Then for any $H \subset G$ as above, we have
$$(V \hat{\otimes}W)^{H-\an} = V^{H-\an}\hat{\otimes}W \quad \textrm{and } (V \hat{\otimes}W)^{\la} = V^{\la}\hat{\otimes}W.$$
\end{prop}
\begin{proof}
We only need to prove the first assertion as the second will follow by taking the inductive limit. By definition, the space $\cal{C}^{\an}(H,V)$ is $\cal{C}^{\an}(H,\Qp)\hat{\otimes}V$. In particular, since the completed tensor product is associative \cite[§2.1 Prop. 6]{BGR}, we get that 
$$\cal{C}^{\an}(H,V \hat{\otimes}W) = \cal{C}^{\an}(H,V) \hat{\otimes}W.$$
Recall that $(V \hat{\otimes}W)^{H-\an} = \cal{C}^{\an}(H,V \hat{\otimes}W)^{\Delta_{1,2}}$. This tells us that
$$(V \hat{\otimes}W)^{H-\an}= (\cal{C}^{\an}(H,V) \hat{\otimes}W)^{\Delta_{1,2}}$$
and since $G$ acts trivially on $W$, this is equal to
$$(\cal{C}^{\an}(H,V))^{\Delta_{1,2}}\hat{\otimes}W=V^{H-\an}\hat{\otimes}W.$$
\end{proof}

The following proposition gives us a sufficient condition for an action on a Banach space to be locally analytic:

\begin{prop}
\label{prop sufficient for locana}
Let $G$ be a $p$-adic Lie group and let $W$ be a $p$-adic Banach representation of $G$. Assume that there exists a compact open subgroup $H$ of $G$ such that, for all $g \in H$, we have
$$||g-1||<p^{-\frac{1}{p-1}}$$
for the operateur norm on $W$. Then the action of $G$ on $W$ is locally analytic.
\end{prop}
\begin{proof}
See \cite[Lemm. 2.14]{BSX18}.
\end{proof}

\subsection{Locally analytic vectors relative to $G_\infty$}
Because of the following result, $\Gal(L/K_{\infty})$ is not necessarily pro-cyclic when $p=2$:
\begin{prop}
\label{careful p=2}
\begin{enumerate}
\item if $K_{\infty} \cap K_{\cycl}=K$, then $\Gal(L/K_{\cycl})$ and $\Gal(L/K_{\infty})$ topologically generate $G_\infty$;
\item if $K_{\infty} \cap K_{\cycl} \supsetneq K$, then necessarily $p=2$, and  $\Gal(L/K_{\cycl})$ and $\Gal(L/K_{\infty})$ topologically generate an open subgroup of $\hat{G}$ of index $2$.
\end{enumerate}
\end{prop}
\begin{proof}
For the first point, see \cite[Lem. 5.1.2]{Liu08} and for the second one, see \cite[Prop. 4.1.5]{Liu10}.
\end{proof}
If $\Gal(L/K_\infty)$ is not pro-cyclic, we let $\Delta \subset \Gal(L/K_\infty)$ be the torsion subgroup, so that $\Gal(L/K_\infty)/\Delta$ is pro-cyclic and we choose $\gamma' \in \Gal(L/K_\infty)$ such that its image in $\Gal(L/K_\infty)/\Delta$ is a topological generator. If  $\Gal(L/K_\infty)$ is pro-cyclic, we choose $\gamma'$ to be a topological generator of $\Gal(L/K_\infty)$.

Let $\tau_n := \tau^{p^n}$ and $\gamma_n:=(\gamma')^{p^n}$. Let $G_n \subset \Gal(L/K)$ be the subgroup topologically generated by $\tau_n$ and $\gamma_n$. It is easy to check that these $G_n$ satisfy the property discussed above lemma \ref{Gn-an subset Gm-an}.

Given a $G_\infty$-representation $W$, we use
$$W^{\tau=1}, \quad W^{\gamma=1}$$
to mean $$ W^{\Gal(L/K_{\cycl})=1}, \quad
W^{\Gal(L/K_{\infty})=1}.$$
And we use
$$
W^{\tau-\la}, \quad W^{\tau_n-\an}, \quad  W^{\gamma-\la} $$
to mean
$$
W^{\Gal(L/K_{\cycl})-\la}, \quad
W^{\Gal(L/(K_{\cycl}(\pi_n)))-\la}, \quad
W^{\Gal(L/K_{\infty})-\la}.  $$

Let
$$ \nabla_\tau  := \frac{\log \tau^{p^n}}{p^n} \text{ for } n \gg0, \quad \nabla_\gamma:=\frac{\log g}{\log_p \chi_p(g)} \text{  for } g \in \Gal(L/K_\infty) \textrm{ close enough to } 1  $$
be the two differential operators acting on $G_\infty$-locally analytic representations.

\begin{rema}
We do not define $\gamma$ as an element of $\Gal(L/K_\infty)$ even though when $\Gal(L/K_\infty)$ is pro-cyclic (and so in particular as soon as $p \neq 2$) we could take $\gamma=\gamma'$. In particular, although the expression ``$\gamma=1$'' might be ambiguous in some cases, we use this notation for brevity.
\end{rema}

Note that if we let $W^{\tau-\la, \gamma=1}:= W^{\tau-\la} \cap W^{\gamma=1}$, then we have $ W^{\tau-\la, \gamma=1} \subset  W^{G_\infty-\la}$ by \cite[Lemm. 3.2.4]{GP18}. We also have $W^{\gamma-\la} \cap W^{\tau=1}\subset  W^{G_\infty-\la}$ since $\Gal(L/K_{\cycl})$ is normal in $\hat{G}$. 

We now recall some results from \cite{GP18} and \cite{Ber14MultiLa} about locally analytic vectors for $G_\infty$ inside some rings of periods. For $n \geq n$, let $r_n = p^{n-1}(p-1)$.

\begin{theo}
\label{theo loc ana basic Kummer case}
Let $I = [r_\ell;r_k]$ or $[0;r_k]$. Then there exists $m_0 \geq 0$, depending only on $k$, such that:
\begin{enumerate}
\item $(\Bt^I_{L})^{\tau_{m+k}-\an, \gamma=1}  \subset \B^I_{\tau,K,m}$ for any $m \geq m_0$;
\item $(\Bt^I_{L})^{\tau-\la, \gamma=1}  = \B^I_{\tau,K,\infty}$;
\item $(\Bt_{\mathrm{rig},L}^{\dagger,r_\ell})^{\tau-\pa, \gamma=1} = \B_{\tau,\mathrm{rig},K,\infty}^{\dagger,r_\ell}$.
\end{enumerate}
\end{theo}
\begin{proof}
This is \cite[Thm. 3.4.4]{GP18}.
\end{proof}

\begin{theo}
\label{theo loc ana cyclo case}
Let $I = [r_\ell;r_k]$ or $[0;r_k]$. Then there exists $m_0 \geq 0$, depending only on $k$, such that:
\begin{enumerate}
\item $(\Bt^I_{L})^{\gamma_{m+k}-\an, \tau=1}  \subset \B^I_{K,m}$ for any $m \geq m_0$;
\item $(\Bt^I_{L})^{\gamma-\la, \tau=1}  = \B^I_{K,\infty}$;
\item $(\Bt_{\mathrm{rig},L}^{\dagger,r_\ell})^{\gamma-\pa, \tau=1} = \B_{\mathrm{rig},K,\infty}^{\dagger,r_\ell}$.
\end{enumerate}
\end{theo}
\begin{proof}
See \cite[Thm. 4.4]{Ber14MultiLa}.
\end{proof}

\begin{theo}
\label{theo loc ana gen Kummer case}
Let $I=[r_\ell;r_k]$ and let $M$ be a finite extension of $K$. Let $M_\infty = M\cdot K_\infty$. Then
\begin{enumerate}
\item $(\Bt_L^I)^{\tau-\la,\Gal(L/M_\infty)=1} = \bigcup_{n \geq 0}\phi^{-n}(\B_{\tau,M}^{p^nI})$;
\item $(\Bt_{\mathrm{rig},L}^{\dagger,r_\ell})^{\tau-\pa,\Gal(L/M_\infty)=1} = \bigcup_{n \geq 0}\phi^{-n}(\B_{\tau,\mathrm{rig},M}^{\dagger,p^nr_\ell})$.
\end{enumerate}
\end{theo}
\begin{proof}
See \cite[Thm. 4.2.9]{GP18}.
\end{proof}

\begin{prop}
\label{invarnabla}
Let $W$ be a $\Qp$-Banach representation of $G_\infty$. Then
$$(W^{\la})^{\nabla_\gamma=0} = \bigcup_{K_\infty \subset_{\mathrm{fin}} M_\infty \subset L}W^{\tau-\la,\Gal(L/M_\infty)=1},$$
where $M_\infty$ runs through the set of all finite extensions of $K_\infty$ inside $L$.
\end{prop}

We now exhibit some ``interesting'' locally analytic vectors for $G_\infty$ inside the rings $\Bt_L^I$. Let $\lambda:= \prod_{n \geq 0}\phi^n(\frac{E([\tilde{\pi}])}{E(0)}) \in \B_{\tau,\mathrm{rig},K}^+$ as in \cite[1.1.1]{KisinFiso}, let $t \in \B_{\mathrm{rig},K}^+$ be the usual $t$ in $p$-adic Hodge theory, and let $b:= \frac{t}{\lambda} \in \At_L^+$, which is exactly $p \cdot \mathfrak{t}$, where $\mathfrak{t}$ is defined in \cite[Ex. 3.2.3]{Liu08}. Note that since $\Bt_L^\dagger$ is a field, there exists some $r(b) \geq 0$ such that $1/b \in \Bt_L^{\dagger,r(b)}$.

\begin{lemm}
If $r_\ell \geq r(b)$, then both $b$ and $1/b$ belong to $(\Bt_L^{[r_\ell,r_k]})^{\la}$.
\end{lemm}
\begin{proof}
See \cite[Lemm. 5.1.1]{GP18}.
\end{proof}

\begin{lemm}
Let $I = [r;s] \subset (0;+\infty)$ such that $r \geq r(b)$. For $n \geq 0$, there exists $b_n \in (\Bt_L^I)^{\la,\nabla_\gamma=0}$ such that $b-b_n \in p^n\At_L^I$.
\end{lemm}
\begin{proof}
This is \cite[Lemm. 5.3.2]{GP18}.
\end{proof}

\section{Kiehl and Tate gluing property}
In this section, we recall results from \cite{DLLZ} that establish Kiehl's gluing property for coherent sheaves on noetherian adic spaces.

We recall the gluing formalism from \cite[Appendix A]{DLLZ}.



Recall the following definition from \cite[Def 1.3.7]{KL01} -

\begin{defi}
\label{GlueDefn}
A gluing diagram is a commuting diagram of ring homomorphisms 

\begin{equation*}
\begin{array}{ccc}

R  & \rightarrow  & R_1 \\
\downarrow &  & \downarrow \\
R_2 & \rightarrow & R_{12} \\

\end{array}
\end{equation*}
such that the $R$-module sequence 

\begin{equation*}
    0 \rightarrow R \rightarrow R_1 \oplus R_2 \rightarrow R_{12} \rightarrow 0
\end{equation*} where the last nontrivial arrow is defined as the difference of the given homomorphisms, is exact.

A gluing datum over this gluing diagram consists of modules $M_1, M_2, M_{12}$ over $R_1, R_2, R_{12}$ respectively, such that there are isomorphisms $\psi_1 : M_1 \otimes R_{12} \cong M_{12}$ and $ \psi_2 : M_2 \otimes R_{12} \cong M_{12} $. Such a datum is said to be finite if the modules are finite over their respective base rings.

\end{defi}

For a gluing datum, we define $M : = ker (\psi_1 - \psi_2 : M_1 \oplus M_2 \rightarrow M_{12})$. There exist natural morphisms $M \rightarrow M_1$ and $M \rightarrow M_2$ of $R$-modules, which induce maps $M\otimes R_1 \rightarrow M_1$ and $M \otimes R_2 \rightarrow M_2$ of $R_1, R_2$-modules respectively.

The following result is \cite[Lem. 1.3.8]{KL01} -

\begin{lemm}
For a finite gluing datum such that $M \otimes R_1 \rightarrow M_1$ is surjective, the following hold true.

\begin{enumerate}
    \item The morphism $\psi_1 - \psi_2 : M_1 \oplus M_2 \rightarrow M_{12} $ is surjective.
    
    \item The morphism $M \otimes R_2 \rightarrow M_2$ is also surjective.
    
    \item There exists a finitely generated $R$-submodule $M_0$ of $M$, such that the morphisms $M_0 \otimes R_1 \rightarrow M_1$, $M_0 \otimes R_2 \rightarrow M_2$ are surjective.
\end{enumerate}
\end{lemm}

The following then is \cite[Lem. A.3]{DLLZ}.

\begin{lemm}
\label{GlueDatumLemma}
In the setting of previous lemma, assume moreover that $R_i$ are noetherian for $i = 1, 2$ and that $R_i \rightarrow R_{12}$ is flat. Suppose that for every finite gluing datum, $M \otimes R_1 \rightarrow M_1$ is surjective. Then for any finite gluing datum, $M$ is a finitely presented $R$-module, and the morphisms $M \otimes R_1 \rightarrow M_1$ and $M \otimes R_2 \rightarrow M_2$ are bijective.
\end{lemm}

We recall the following definition from \cite{DLLZ}.

\begin{defi}
We call a homomorphism of Huber rings $f : R \rightarrow S$ strict adic if, for one, and hence for every choice of an ideal of definition $I \subset R$, its image $f(I)$ is an ideal of definition for $S$. It follows that a strict adic homomorphism is an adic homomorphism.
\end{defi}

For a strict adic homomorphism, we have the following gluing lemma. (See \cite[Lem. 2.7.2]{KL01}, \cite[Lem. A.6]{DLLZ}.)

\begin{lemm}
\label{StrictAdic}
Let $R_1 \rightarrow S$ and $R_2 \rightarrow S$ be adic homomorphisms of complete Huber rings such that their direct sum $\psi : R_1 \oplus R_2 \rightarrow S$ is a strict adic homomorphism. Then, for any ideal of definition $I_S \subset S$, there exists some $l > 0$ such that, every $U \in GL_n(S)$ satisfying $U - 1 \in M_n(I_{S}^{l})$ can be written as $\psi(U_1) \psi(U_2)$ for $U_i \in GL_n(R_i)$.
\end{lemm}

\begin{proof}
We briefly recall the proof. By hypotheses, for any ideals of definitions $I_1 \subset R_1$ and $I_2 \subset R_2$, we have $I^{'}_{S} : = \psi_{I_1 \oplus I_2} \subset S$ as an ideal of definition. We can choose $l > 0 $ such that $I^{l}_{S} \subset I^{'}_{S}$ since both are ideals of definition. Then it suffices to prove that any $U \in GL_n(S)$ satisfying $U - 1 \in M_n(I_{S}^{'})$ can be written as $\psi(U_1) \psi(U_2)$ for $U_i \in GL_n(R_i)$. Given $U \in GL_n(S) $ with $ V : = U - 1 M_n(I^{'m}_{S})$ for some $m > 0$, we know by assumption that $V$ arises from a pair $(X, Y) \in M_{n}(I_{1}^{m}) \times M_{n}(I_{2}^{m})$. Then it follows that $U^{'} := \psi(1-X)U\psi(1-Y)$ satisfies $U^{'} - 1 \in M_{n}(I^{'2m}_{S})$ and we conclude by iterating this procedure. 
\end{proof}

The following key lemma then forms the heart of the gluing argument. (\cite[Lem. 2.7.4]{KL01}, \cite[Lem. A.7]{DLLZ}.)

\begin{lemm}
\label{KeyGluing}
In the context of definition \ref{GlueDefn} and the definition of $M$, assume furthermore that :

\begin{enumerate}
    \item The Huber rings $R_1, R_2, R_{12}$ are complete;
    \item $R_1 \oplus R_2 \rightarrow R_{12} $ is a strict adic homomorphism; and
    \item The map $R_2 \rightarrow R_{12}$ has dense image.
\end{enumerate}
Then for $i = 1, 2$, the natural map $M \oplus R_i \rightarrow M_i $ is surjective.
\end{lemm}

\begin{proof}
Choose sets of generators $\{ m_{1, 1}, m_{2, 1}, \ldots , m_{n, 1} \}$ and $\{ m_{1, 2}, m_{2, 2}, \ldots , m_{n, 2} \}$ of $M_1$ and $M_2$ respectively of the same cardinality. Then there exist matrices $A, B \in M_n(R_{12})$ such that $$\psi_{2}(m_{j, 1}) = \sum_{i} A_{ij} \psi_{1}(m_{i, 2})$$ and $$\psi_{1}(m_{j, 2}) = \sum_{i} B_{ij} \psi_{2}(m_{i, 1})$$ for every $j$.

Since $R_2 \rightarrow R_{12} $ has dense image, by Lemma \ref{StrictAdic}, there exists matrix $B^{'} \in M_{n}(R_{2})$ such that $1 + A(B^{'} - B) = C_1 C^{-1}_{2}$ for some $C_{i} \in R_{i}$. For each $j = 1, 2, \ldots, n$, let $x_j : = (x_{j, 1}, x_{j, 2}) = (\sum_{i} (C_1)_{ij} m_{i, 1}, \sum_{i} (B^{'}C_2)_{ij} m_{i, 2}) \in M_1 \times M_2$. Then \begin{equation*}
    \psi_{1}(x_{j,1}) - \psi_{2}(x_{j, 2}) = \sum_{i} (C_1 - AB'C_2)_{ij} \psi_{1}(m_{i, 1})
\end{equation*} and hence \begin{equation*}
    \psi_{1}(x_{j,1}) - \psi_{2}(x_{j, 2}) = \sum_{i} ((1-AB)C_2)_{ij} \psi_{1}(m_{i, 1}) = 0 
\end{equation*} by definition of $A $ and $B$. Thus, $x_{j} \in M$. But $C_i$ are invertible and hence $\{ x_{j, i} \}_{j = 1}^{n} $, for $i = 1, 2$, generates $M_i$ over $R_i$, thus giving the claim.
\end{proof}

Finally, we come to the theorem we need for our application.

\begin{theo}
\label{GlueThm}
Let $X = \mathrm{Spa } (R, R^{+})$ be a noetherian affinoid adic space. Then the categories of coherent sheaves on $X$ and finitely generated $R$-modules are equivalent to each other via the global sections functor.
\end{theo}

\begin{proof}
It suffices to check the Kiehl gluing property on simple Laurent coverings $\mathrm{Spa} (R_{i}, R^{+}_{i}), i = 1,2$ by \cite[Lem. 2.4.20]{KL01}. For any such covering, define $\mathrm{Spa}(R_{12}, R^{+}_{12}) : = \mathrm{Spa} (R_{1}, R^{+}_{1}) \times \mathrm{Spa} (R_{2}, R^{+}_{2})$. (Recall that coproducts exist in the category of adic spaces.) By the Noetherian hypothesis and \cite[Thm. 2.5]{Huber}, $R, R_i, R_{12}$ form a gluing diagram. Further, $R_i \rightarrow R_{12}$ is flat with dense image for $i = 1, 2$. Thus, we can conclude by applying Lemmas \ref{GlueDatumLemma} and \ref{KeyGluing}.
\end{proof}

\section{Families of representations and $(\phi,\Gamma)$-modules}
\subsection{Families of representations}
We let $S$ be a $\Qp$-Banach algebra, and we let $\cal{X}$ be the set of maximal ideals of $S$. As in \cite{BC08}, we think of elements of $\cal{X}$ as points and we write $\mathfrak{m}_x$ for the maximal ideal of $S$ corresponding to a point $x \in \cal{X}$. For $f \in S$, we let $f(x)$ denote the image of $f$ in $E_x = S/\mathfrak{m}_x$. 

Instead of working with norms, we work with ``valuations'' on $S$, such that for any $f,g \in S$, we have $\mathrm{val}_S(fg) \geq \mathrm{val}_S(f)+\mathrm{val}_S(g)$. 

Following \cite[\S 2]{BC08}, we say that $S$ is an algebra of coefficients if $S$ satisfies the following conditions:
\begin{enumerate}
\item $S \supset \Qp$ and the restriction of $\mathrm{val}_S$ to $\Qp$ is the $p$-adic valuation $v_p$;
\item for any $x \in \cal{X}, E_x$ is a finite extension of $\Qp$;
\item the Jacobson radical $\mathrm{rad}(S)$ is zero.
\end{enumerate}

Let $S$ be a $\Qp$-Banach algebra. A family of $p$-adic representations of $\G_K$ is an $S$-module $V$ free of finite rank $d$, endowed with a continuous linear action of $\G_K$. Under the assumption that there exists a free $\O_S$-module (where $\O_S$ is the ring of integers of $S$ for $\mathrm{val}_S$) $T$ of rank $d$ such that $V=S \otimes_{\O_S}T$, Berger and Colmez show in \cite{BC08} how to attach to such a family of representations a family of $(\phi,\Gamma)$-modules over $S \hat{\otimes}\B_K^\dagger$, using what are called Sen-Tate conditions. They also use a result of étale descent which we will also need:

Let $B$ be a $\Qp$-Banach algebra endowed with a continuous action of a finite group $G$. Let $B^\natural$ denote the ring $B$ endowed with the trivial $G$-action, and assume that:
\begin{enumerate}
\item the $B^G$-module $B$ is finite free and faithfully flat;
\item we have $B^\natural \otimes_{B^G}B \simeq \oplus_{g \in G}B^\natural \cdot e_g$ (where $e_g^2=e_g, e_ge_h=0$ if $g \neq h$ and $g(e_h)=e_{gh}$). 
\end{enumerate}

\begin{prop}
\label{prop classical etale descent}
If $S$ is a Banach algebra (on which $G$ acts trivially), and if $M$ is an $S\hat{\otimes}B$-module locally free of finite type endowed with a semilinear action of $G$ then:
\begin{enumerate}
\item $M^G$ is an $S\hat{\otimes}B^G$-module locally free of finite type;
\item the map $(S\hat{\otimes}B)\otimes_{S\hat{\otimes}B^G}M^G \rightarrow M$ is an isomorphism.
\end{enumerate}
\end{prop}
\begin{proof}
This is \cite[Prop. 2.2.1]{BC08}.
\end{proof}

\subsection{Tate-Sen formalism for Huber rings}
Here we formulate Tate-Sen formalism for Huber rings. This was developed by the first author in joint work with Ruochuan Liu (\cite{KarnatakiLiu21}). In \cite{BC08} this is done for $\Q_p$-Banach algebras but the generalization to Huber rings is straightforward.

Recall that $A$ is called Huber ring if there exists an open adic subring $A_0 \subset A$ (called ring of definition of $A$) with finitely generated ideal of definition $I$. We recall the notion of boundedness for Huber rings.

\begin{defi}
Let $A$ be a huber ring. A subset $\Sigma \subset A$  is bounded if for every open neighbourhood $U$ of $0$ in $A$, there exists an open neighbourhood $V$ of $0$ in $A$ such that $$ V . \Sigma \subset U.$$ 
\end{defi}

We note that any ring of definition $A_0 \subset A$ has to be bounded. Conversely, any open bounded subring of $A$ is a ring of definition for $A$. In the case of $k$ an archimedian field and $A$ being a reduced affinoid algebra over $k$, the set $A^{0}$ of power-bounded elements is the (closed) unit ball under the sup-norm. 

Now let $A$ be a Huber ring, and $\tilde{S}$ be an $A^0$-algebra. We denote by $I_{\tilde{S}}$ the ideal of definition for $\tilde{S}$. We state the generalised Tate-Sen formalism in this setting below. 

\begin{defi}
Let $G$ be a group acting on an adic ring $R$. We say $G$ acts on $R$ strict adically if, for each $g \in G$, the action of $g$ on $R$ gives a strict adic homomorphism $R \rightarrow R$.
\end{defi}

Assume that $\tilde{S}$ has an action of a profinite group $G_0$. We assume that it acts on $\tilde{S}$ strict adically. As before, we also fix a character $\chi : G_0 \rightarrow \mathbb{Z}^{\times}_{p}$ with open image and set $H_0 : = \mathrm{ker} \chi$. For any open subgroup $G \subset G_0$, we define $H := G \cap H_0$. We set $G_H$ to be the normaliser of $H$ in $G_0$, and we define $\tilde{\Gamma}_H := G_H / H$.

Then the Tate-Sen conditions are as follows.

\begin{enumerate}[TS1]
    \item There exists an integer $l_1 > 0 $, such that for any open subgroups $H_1 \subset H_2$ of $H_0$, there exists an $\alpha \in \tilde{S}^{H_1}$ such that $\alpha . I_{\tilde{S}}^{l_1} \subset \tilde{S}_0$ and $\sum_{\tau \in H_2 / H_1} \tau(\alpha) = 1$. 
    
    \item For each open subgroup $H$ of $H_0$, there exists an increasing sequence $(S_{H, n})_{n \ge 0}$ of closed sub-$S^0$-algebras of $\tilde{S}^{H}$, and an integer $n(H) \ge 0$ such that for each $n \ge n(H)$, there is an $S^0$-linear map $R_{H, n} : \tilde{S}^{H} \rightarrow S_{H, n} $. There is also an integer $l_2 $ independent of $H$, such that the following properties are satisfied by this collection of objects. 
    
    \begin{enumerate}[(a)]
            \item For $H_1 \subset H_2$, we have $S_{H_2, n} \subset S_{H_1, n}$ and the restriction of $R_{H_1, n}$ to $\tilde{S}^{H_2}$ coincides with $R_{H_2, n}$.
            
            \item $R_{H, n}$ is $S_{H, n}$-linear and $R_{H, n}(x) = x$ if $x \in S_{H, n}$.
            
            \item $g(S_{H, n}) = S_{gHg^{-1}, n}$ and $g(R_{H, n}(x)) = R_{gHg^{-1}, n}(gx)$ for all $g \in G_0$; in particular, $R_{H, n}$ commutes with the action of $\tilde{\Gamma}_H$.
            
           \item If $n \ge n(H)$, and $x \in I_{\tilde{S}}^{m}$, then $R_{H, n}(x) \in I_{\tilde{S}}^{m - l_2}$.
            
            \item If $x \in \tilde{S}^{H}$, then $\mathrm{lim}_{n \rightarrow \infty} R_{H, n}(x) = x$.
        \end{enumerate}
    
    \item There exists an integer $l_3 > 0$, and for each open subgroup $G \subset G_0$, an integer $n(G) \ge n_1(H)$, where $H = G \cap G_0$, such that if $n \ge n(G)$ and $\gamma \in \tilde{\Gamma}_{H} $ satisfies $n(\gamma) < n$, then $\gamma - 1$ is invertible on $X_{H, n} : = (1 - R_{H. n}) (S_{H, n})$, and if $x \in I_{\tilde{S}}^m$, then $ (\gamma - 1)^{-1}(x) \in I_{\tilde{S}}^{m- l_3}$ for all $x \in X_{H, n}$.
    
\end{enumerate}

For an ensemble of objects satisfying these axioms, we prove that an analogue of the theorem of Berger and Colmez holds.

\begin{theo}[Existence of $(\varphi, \Gamma)$-modules over an appropriate Robba ring]
\label{thm:PhiGammaDescent}

Let $A$ be a Huber ring and let $\tilde{S}$ be an $A^0$-algebra satisfying $(TS1), (TS2),$ and $(TS3)$ as above. Let $T$ be an $A^0$-representation of dimension $d$ of $G_0$, $V = A \otimes_{A^0} T$, and $k$ be an integer such that $p^k \in I_{\tilde{S}}^{l_1 + 2l_2 + 2l_3}$. Let $G$ be the subgroup of $G_0$ acting trivially on $T/ p^kT$, let $H = G \cup H_0$ and let $n \ge n(G)$. Then, $\tilde{S}^{0} \otimes_{A^0} T$ contains a unique $S_{H. n}^0$-submodule $D^{0}_{H, n}(T)$, which is free of rank $d$ and satisfies the following properties -

\begin{enumerate}
    \item $D^{0}_{H, n}(T)$ is fixed by $H$ and stable unde the action of $G_0$.
    
    \item The natural map $$D^{0}_{H, n}(T) \otimes_{S_{H, n}^{0}} \tilde{S}^0 \rightarrow \tilde{S} \otimes_{A^0} T$$ is an isomorphism.
    
    \item There is a basis of $D^{0}_{H, n}(T)$ over $S_{H, n}^{0}$ that is $l-3$-fixed by $G/H$. That is, for any $\gamma \in G/H$, the matrix $W$, by which $\gamma$ acts in this basis, belongs to $M_d(I_{S_{H, n}}^{l_3})$.
\end{enumerate}
\end{theo}

We first prove a number of lemmas needed for this proof, deferring the proof to the end of this section.

\begin{lemm}
Let $H$ be an open subgroup of $H_0$. If $a > l_1$ is an integer, and $k \in \mathbb{N}$ and if $\tau \rightarrow U_{\tau}$ is a continuous cocycle of $H$ valued in $GL_d(\tilde{S})$ satisfying $U_{\tau} - 1 \in p^k M_d(\tilde{S})$, and $U_{\tau} - 1 \in M_{d}(I_{\tilde{S}}^a)$, then for all such $\tau \in H$, there exists a matrix $M \in GL_d(\tilde{S})$, satisfying $M - 1 \in p^k M_d(\tilde{S})$, and $M - 1 \in M_d(I_{\tilde{S}}^{a - l_1})$ such that the cocycle $\tau \rightarrow M^{-1} U_{\tau} \tau(M)$ satisfies $(M^{-1} U_{\tau} \tau(M) - 1 ) \in I_{\tilde{S}}^{a+1}$.
\end{lemm}

\begin{proof}
Let $H_1$ be an open subgroup of $H$ such that $U_{\tau} - 1 \in I_{\tilde{S}}^{a + 1 + l_1}$ if $\tau \in H_1$. Let $\alpha \in \tilde{S}^{H_1}$ such that $\sum_{\tau \in H/H_1} \tau(\alpha) = 1$ and $\alpha . I_{\tilde{S}}^{l_1} \subset \tilde{S}_0$. If $Q$ is a system of representatives for $H/H_1$, we define $$ M_Q = \sum_{\sigma \in Q} \sigma(\alpha)U_{\sigma}.$$ We have $M_Q - 1 = \sum_{\sigma \in Q} \sigma(\alpha)(U_{\sigma} - 1)$. This implies that $M_Q - 1 \in M_{d}(I_{\tilde{S}}^{a - l_1})$. Moreover, $M_{Q}^{-1} = \sum_{n=0}^{\infty} (1 - M_Q)^n$, since the sum on right hand side converges, and $M_{Q}^{-1} \in M_d(I_{\tilde{S}}^{m})$ for some $m \ge 0$ and $M_Q \in GL_d(\tilde{S})$.

If $\tau \in H_1$, then by the cocycle condition we get $U_{\tau \sigma} - U_{\sigma} = U_{\sigma}(\sigma(U_{\tau}) - 1) $. Let $Q'$ be another set of representatives for $H / H_1$. Then, for any $\sigma' \in Q'$ there exists a $\tau \in H_1$ and $\sigma \in Q$ such that $\sigma' = \sigma \tau$. Thus, we get $$ M_Q - M_{Q'} = \sum_{\sigma \in S} \sigma(\alpha) (U_{\sigma} - U_{\sigma \tau}) = \sum_{\sigma \in S} \sigma(\alpha) U_{\sigma} (1 - \sigma(U_{\tau})).$$ Thus, $$M_Q - M_{Q'} \in M_{d}(I_{\tilde{S}}^{a+1}).$$ For any $\tau \in H_0$, $$U_{\tau} \tau(M_Q) = \sum_{\sigma \in Q} \tau \sigma(\alpha) U_{\tau} \tau(U_{\sigma}) = M_{\tau Q}.$$ Then, $$M_{Q}^{-1} U_{\tau} \tau(M_Q) = 1 + M_{Q}^{-1} (M_{\tau Q} - M_Q) $$ with $M_{Q}^{-1} (M_{\tau Q} - M_Q) \in M_{d}(I_{\tilde{S}}^{a+1})$. Setting $M = M_Q$, we get the result.
\end{proof}

\begin{coro}
\label{cor:descent}

Under the same hypotheses as the above lemma, there exists a matrix $M \in GL_{d}(\tilde{S})$ such that $M \in M_{d}(I_{\tilde{S}}^{a - l_1})$ and $$ M^{-1}U_{\sigma} \sigma(M) = 1 $$ for all $\sigma \in H_0$.
\end{coro}

\begin{proof}
Repeat the lemma for $(a \rightarrow a+1 \rightarrow a+2 \rightarrow \cdots)$ and take limits of the matrices you get from the lemma. 
\end{proof}

\begin{lemm}
Let $\delta > 0$. Let $a, b \in \mathbb{R}$ such that $a \ge l_2 + l_3 + \delta$ and $ b \ge \mathrm{Sup}(a + l_2, 2l_2+2l_3+\delta)$. Let $H$ be an open subgroup of $H_0$, $n \ge n(H)$, $\gamma \in \tilde{\Gamma}_{H}$ satisfying $n(\gamma) \le n$ and let $ U = 1 + p^kU_1 + p^kU_2$, with - $$ U_1 \in M_d(I_{S_{H, n}}^{a - r}), U_2 \in M_d(I_{\tilde{S}^{H}}^{b - r})$$ where $r := \mathrm{max} (n) : p^k \in I_{\tilde{S}}^n$.

Then, there exists a matrix $M \in 1 + p^kM_d(\tilde{S}^{H})$ such that $M - 1 \in I_{\tilde{S}}^{b - l_2 - l_3}$ such that $M^{-1} U \gamma(M) = 1 + p^kV_1 + p^kV_2$ with -  $$ V_1 \in M_d(I_{S_{H, n}}^{a - r}), V_2 \in M_d(I_{\tilde{S}^{H}}^{b - r + \delta}).$$
\end{lemm}

\begin{proof}
By the conditions $(TS2)$ and $(TS3)$, we can write $U_2$ in the form $U_2 = R_{H, n}(U_2) + (1 - \gamma)(V)$ with $R_{H, n}(p^kU_2) \in M_d(I_{\tilde{S}}^{b - l_2})$ and $p^kV \in M_d(I_{\tilde{S}}^{b - l_2 - l_3})$. 

Thus, $$ (1+p^kV)^{-1} U \gamma(1 + p^kV) = (1 - p^kV + p^{2k}V^2 - \ldots) (1 +  p^kU_1 + p^kU_2) (1 + p^k \gamma(V)).$$ This gives $$ (1+p^kV)^{-1} U \gamma(1 + p^kV) = 1 + p^kU_1 + (\gamma - 1)V + p^kU_2 + (\text{terms of degree } \ge 2).$$

Let $V_1 = p^kU_1 + p^kR_{H, n}(U_2)$ and let the terms of degree $ \ge 2$ be denoted by $W$. Then we see that $M = (1 + p^kV) $ and $V_2 = W$ gives us the result.
\end{proof}

\begin{coro}
\label{cor:decompletion}
Under the same hypotheses as the above lemma, there exists a matrix $M \in GL_d(\tilde{S}^{H})$ with further $M -1 \in M_d(I_{\tilde{S}}^{b - l_2 - l_3})$ such that $ M^{-1} U \gamma(M) \in GL_{d}(S_{H, n}) $.
\end{coro}

\begin{proof}
Repeat the above construction for $(b \rightarrow b+\delta \rightarrow b + 2\delta \rightarrow \ldots) $ (Take $\delta = 1$ in fact). Then take the limit.
\end{proof}

\begin{lemm}
\label{lem:translate}

Let $H$ be an open subgroup of $H_0$, $n \ge n(H)$. Let $\gamma \in \tilde{\Gamma}_{H}$ satisfying $n \le n(\gamma)$ and let $B \in GL_d(\tilde{S}^{H})$. If there exist $V_1, V_2 \in GL_d(S_{H, n})$ with $V_1 - 1, V_2 - 1 \in M_d(I_{\tilde{S}}^{l_3})$ such that $\gamma(B) = V_1 B V_2 $, then $B \in GL_d(S_{H, n})$.
\end{lemm}

\begin{proof}
If $C = B - R_{H, n}(B)$, then $\gamma(C) = V_1CV_2$, since the map $R_{H,n}$ is $S_{H,n}$-linear and commutes with the action of $\gamma$. We have to prove $C = 0$. We have - $$ \gamma(C) - C = V_1CV_2 - C = (V_1 - 1)CV_2 + V_1C(V_2 - 1) + (V_1 - 1)C(V_2 -1).$$ Hence, if $C \in M_d(I_{\tilde{S}}^{m})$, we have $ \gamma(C) - C \in M_d(I_{\tilde{S}}^{m + l_3+1})$ which by $(TS3)$ implies that $C = 0$.
\end{proof}

Finally, we come to the proposition that connects all the lemmas together.

\begin{prop}
\label{prop:descend}

Let $\tilde{S}$ be an $f$-adic ring satisfying the axioms $(TS1), (TS2), (TS3)$ for $f$-adic rings. Let $\sigma \rightarrow U_{\sigma}$ be a continuous $1$-cocycle for $G_0$ taking values in $GL_d(\tilde{S})$. If $G$ is a distinguished open subgroup of $G_0$ such that $U_{\sigma} - 1 \in p^kM_{d}(\tilde{S})$, and in fact $U_{\sigma} - 1 \in M_{d}(I_{\tilde{S}}^{l_1 + 2l_2 + 2l_3})$ for all $\sigma \in G$ and if $H = H_0 \cap G$, then there exists a matrix $M$ such that $M \in 1 + p^kM_{d}(\tilde{S})$ satisfying $M - 1 \in M_{d}(I_{\tilde{S}}^{l_2 + l_3})$ such that the $1$-cocycle $\sigma \rightarrow V_{\sigma} := M^{-1}U_{\sigma}\sigma(M)$ is trivial over $H$ and takes values in $S_{H,n(G)}$.
\end{prop}

\begin{proof}
Corollary \ref{cor:descent} gives a matrix $M_1 \in 1 + p^kM_d(\tilde{S})$ with $M - 1 \in M_d(I_{\tilde{S}}^{2l_2 + 2l_3})$ such that the $1$-cocycle $\sigma \rightarrow U'_{\tau} := M_1^{-1}U_{\tau} \tau(M_1)$ is trivial on $H$ and thus by inflation provides a $1$-cocycle for the group $\tilde{\Gamma}_H$ taking values in $\tilde{S}^{H}$. (Since $G$ is distinguished in $G_0$, this implies that $G_H = G_0$.)

Let $\gamma \in \tilde{\Gamma}_H$ with $n(\gamma) = n(G)$. In particular, $\gamma$ is in the image of $G$ and $U_{\gamma} - 1 \in p^kM_{d}(\tilde{S}^H)$ with further $U'_{\gamma} - 1 \in M_{d}(I_{\tilde{S}}^{2l_2 + 2l_3})$. By corollary \ref{cor:decompletion}, we get a matrix $M_2 \in 1 + p^kM_d(\tilde{S}^H)$ with $M_2 - 1 \in M_d(I_{\tilde{S}}^{l_2 + l_3})$ such that $M_2^{-1} U'_{\gamma} \gamma(M_2) \in GL_{d}(S_{H, n(G)})$.

Then, letting $M = M_1M_2$, we have $M \in 1 + p^kM_d(\tilde{S})$ and in fact, $M - 1 \in M_d(I_{\tilde{S}}^{l_2 + l_3})$, and the cocyle $ \tau \rightarrow M^{-1} U_{\tau} \tau(M)$ is trivial over $H$ and takes values in $GL_d(\tilde{S}^H)$. In fact, the matrix $V_{\gamma} \in GL_d(S_{H, n(G)})$ and $V_{\gamma} - 1 \in M_{d}(I_{\tilde{S}}^{l_2 + l_3})$. 

It remains to prove that $V_{\tau} \in GL_d(S_{H, n(G)})$ for all $\tau \in G_0$. To this end, if $\tau \in G_0$, we have the relation $\tau \gamma = \gamma \tau$ in $G_0 / H$ and the cocycle condition gives the relation - $$ V_{\tau} \tau(V_{\gamma}) = V_{\gamma} \gamma(V_{\tau}).$$ Then we apply lemma \ref{lem:translate} with $B = V_{\tau}, V_1 = V_{\gamma}^{-1}$ and $V_2 = \tau(V_{\gamma})$ to deduce the fact that $V_{\tau}$ takes values in $GL_d(S_{H, n(G)})$. This finishes the proof.
\end{proof}

We use these results to supply a proof of Theorem \ref{thm:PhiGammaDescent} below.

\begin{proof}[Proof of Theorem \ref{thm:PhiGammaDescent}]
Let $v_1, \ldots, v_d$ be a basis for $T$ over $A^0$ and let $U_{\sigma} = (u^{\sigma}_{i,j})$ be the matrix of vectors $\sigma(v_1), \ldots, \sigma(v_d)$ over the basis $v_1, \ldots, v_d$. Then $\sigma \rightarrow U_{\sigma}$ is a continuous $1$-cocycle taking values in $GL_d(A^0) \subset GL_d(\tilde{S}^0)$.

From the hypotheses, we have $U_{\sigma} \in 1 + p^kM_d(A^0)$ if $\sigma$ is in $G$. By proposition \ref{prop:descend}, we get a matrix $M \in GL_d(\tilde{S})$ satisfying $ M - 1 \in M_d(I_{\tilde{S}}^{l_2 + l_3})$ (and thus $M \in GL_d(\tilde{S}^0)$) such that the cocycle $\sigma \rightarrow V_{\sigma} := M^{-1} U_{\sigma} \sigma(M)$ is trivial over $H$, and takes values in $GL_d(S_{H, n(G)}) \cap GL_d(\tilde{S}^0) = GL_d(S_{H, n(G)}^0)$. If $M = (m_{i,j})$, and if $e_k = \sum_{j=1}^{d} m_{j, k}v_j$, we have $$ \sigma(e_k) = \sum_{j=1}^{d} \sigma(m_{j, k}) \sigma(v_j) = \sum_{i=1}^{d} \left( \sum_{j=1}^{d} u^{\sigma}_{i,j} \sigma(m_{j, k}) \right) v_i = e_k. $$ If $\sigma \in H$, this gives the fact that $e_1, \ldots, e_d$ is a basis for $\tilde{S}^0 \otimes_{A^0} T$ over $\tilde{S}^0$ that is fixed by $H$.

Now, if $\gamma \in G/H$, the matrix $W$ of $\gamma$ in the basis $e_1, \ldots, e_d$ is of the form $M^{-1}U_{\sigma} \sigma(M)$, where $\sigma \in G$ is a lift of $\gamma$, and $W - 1 \in M_d(I_{\tilde{S}}^{l_2 + l_3})$. Thus we deduce that the sub-$S_{H, n(G)}^{0}$-module generated by $e_1, \ldots, e_d$ satisfies the required properties, and thus we get the existence of such a module.

It remains to show the uniqueness. Fix a $\gamma \in \tilde{\Gamma}_{H}$ satisfying $n(\gamma) = n(G)$. Let $e_1, \ldots, e_d$ and $e'_{1}, \ldots, e'_{d}$ are two bases of $\tilde{S}^0 \otimes_{A^0} T$ over $\tilde{S}^0$ fixed by $H$, such that the matrices $W$ and $W'$ of $\gamma$ in these bases are in $GL_d(\tilde{S}^0)$, with $n \ge n(G)$, satisfying $W - 1, W' - 1 \in M_d(I_{\tilde{S}}^{l_3})$. Then, let $B$ be the matrix of the vectors $e'_{j}$ in the basis $e_1, \ldots, e_d$. Then, $B$ is fixed by $H$, and we have $W' = B^{-1} W \gamma(W)$. Then, by lemma \ref{lem:translate}, we deduce that $B$ takes values in $S_{H, n}$, and hence in $S_{H, n}^0$. This implies that the two $S_{H, n}^0$-modules generated by $e_1, \ldots, e_d$ and $e'_{1}, \ldots, e'_{d}$ are the same. This finishes the proof.

\end{proof}

\begin{rema}
Assume the hypotheses of Theorem \ref{thm:PhiGammaDescent}. If we define $ D_{H, n}(T) := S_{H, n} \otimes_{S_{H, n}^0} D^{0}_{H, n}(T)$, then $D_{H, n}(T)$ is a free module of rank $d$ over $S_{H, n}$. It is the unique sub-$S_{H, n}$-module of $\tilde{S} \otimes_{A^0} T$ satisfying the following properties :

\begin{enumerate}
    \item $D_{H, n}(T)$ is fixed by $H$ and is stable under $G_0$.
    
    \item The natural map $ \tilde{S} \otimes_{S_{H, n}} D_{H, n}(T) \rightarrow \tilde{S} \otimes_{A^0} T$ is an isomorphism.
    
    \item $D_{H, n}(T)$ has a basis over $S_{H, n}$ that is $l_3$-fixed by $G/H$. 
\end{enumerate}
The proof is exactly the same as that of Theorem \ref{thm:PhiGammaDescent}.
\end{rema}

\subsection{Overconvergent families of $(\phi,\Gamma)$-modules}
The Sen-Tate conditions applied to the ring $\At^{\dagger,1}$ allow Berger and Colmez to attach a family of $(\phi,\Gamma)$-modules to a family of representations:

Let $S$ be a $\Qp$-Banach algebra, let $V$ be an $S$-representation of $\G_K$, let $T$ be an $\O_S$-lattice of $V$ stable under the action of $\G_K$, and let $M$ be a finite Galois extension of $K$ such that $\G_M$ acts trivially on $T/12pT$. Let $n(M)$ be as defined in \cite[§4]{BC08} and let $r(V) = \max ((p-1)p^{n(M)-1},r'(M))$. Up to increasing $r(V)$, we make sure that there exists $n(V)$ such that $p^{n(V)-1}(p-1) = r(V)$. We also let, as in \cite{BC08}, $c_1$, $c_2$ and $c_3$ be constants such that $c_1 > 0$, $c_2 > 0$, $c_3 > 1/(p-1)$ and such that $c_1+2c_2+2c_3 < v_p(12p)$.

\begin{prop}
\label{prop action of gamma loc ana on phigamma}
If $V$ is an $S$-representation of $\G_K$ of dimension $d$, and if $n \geq n(M)$, then $(\O_S \hat{\otimes}\At^{\dagger,1})\otimes_{\O_S}T$ contains a unique sub-$\O_S\hat{\otimes}\A_{M,n}^{\dagger,1}$-module $D_{M,n}^{\dagger,1}(T)$, free of rank $d$, fixed by $H_M$, stable under $\G_K$ and having an almost $\Gamma_M$-invariant basis such that:
$$(\O_S \hat{\otimes}\At^{\dagger,1})\otimes_{\O_S\hat{\otimes}\A_{M,n}^{\dagger,1}}D_{M,n}^{\dagger,1}(T) \simeq (\O_S\hat{\otimes}\At^{\dagger,1})\otimes_{\O_S}T.$$
Moreover, there exists a basis of $D_{M,n}^{\dagger,1}(T)$ in which, if $\gamma \in \Gamma_M$, then the matrix of $\gamma$ in this basis satisfies $V(W_\gamma-1,[1,+\infty]) > c_3$.
\end{prop}
\begin{proof}
The first part of the proposition is \cite[Prop. 4.2.8]{BC08}. The part on the action of $\Gamma_M$ comes from \cite[Prop. 4.2.1]{BC08} and the Tate-Sen conditions.
\end{proof}

If $V$ is an $S$-representation of $\G_K$ of dimension $d$ and if $r \geq r(V)$, then Berger and Colmez define:
$$D_K^{\dagger,r}(V) = (S\hat{\otimes}\B_M^{\dagger,r} \otimes_{S\hat{\otimes}\B_M^{\dagger,r(V)}}\phi^{n(V)}(D_{M,n(V)}^{\dagger,1}(V)))^{H_K}.$$

Berger and Colmez then prove the following:
\begin{theo}
\label{theo overconvergence phigamma}
If $V$ is an $S$-representation of $\G_K$ of dimension $d$ and if $r \geq r(V)$, then:
\begin{enumerate}
\item $D_K^{\dagger,r}(V)$ is a locally free $S\hat{\otimes}\B_K^{\dagger,r}$-module of rank $d$;
\item the map $(S\hat{\otimes}\Bt^{\dagger,r}) \otimes_{S\hat{\otimes}\B_K^{\dagger,r}}D_K^{\dagger,r} \rightarrow (S\hat{\otimes}\Bt^{\dagger,r})\otimes_S V$ is an isomorphism.
\end{enumerate}
\end{theo}
\begin{proof}
See \cite[Thm. 4.2.9]{BC08}.
\end{proof}

\section{Analytic families of $(\phi,\tau)$-modules over the Robba ring}
\label{section ana families}
In this section, we explain how to attach to a family of $(\phi,\Gamma)$-modules over $S \hat{\otimes}\B_K^\dagger$ (satisfying some additional conditions) a family of $(\phi,\tau)$-modules over $(S \hat{\otimes}\B_{\tau,\rig,K}^\dagger,S \hat{\otimes}\Bt_{\rig,L}^\dagger)$. In particular, using the family of $(\phi,\Gamma)$-modules given by theorem \ref{theo overconvergence phigamma} attached to a family of representations $V$ will give us a family of $(\phi,\tau)$-modules over the Robba ring $S \hat{\otimes}\B_{\tau,\rig,K}^\dagger$ which is canonically attached to $V$. 
\subsection{Overconvergent $(\phi,\Gamma)$-modules and locally analytic vectors}
Let $S$ be a $\Qp$-Banach algebra and let $V$ be an $S$-representation of $\G_K$ of dimension $d$. For $0 \leq r \leq s$, we let
$$\tilde{D}_L^{[r;s]}(V) = ((S\hat{\otimes}\Bt^{[r;s]})\otimes_S V)^{\G_L} \quad \textrm{and } \tilde{D}_{\mathrm{rig},L}^{\dagger,r}(V) = ((S\hat{\otimes}\Bt_{\mathrm{rig}}^{\dagger,r})\otimes_S V)^{\G_L}.$$
These two spaces are topological representations of $G_\infty$. By theorem \ref{theo overconvergence phigamma}, we have an other description of $\tilde{D}_L^{[r;s]}(V)$ and $\tilde{D}_{\mathrm{rig},L}^{\dagger,r}(V)$ for $r \geq s(V)$:
\begin{itemize}
\item $\tilde{D}_L^{[r;s]}(V) = (S\hat{\otimes}\Bt_L^{[r;s]}) \otimes_{S\hat{\otimes}\B_K^{^\dagger,r}}D_K^{\dagger,r}(V)$;
\item $\tilde{D}_{\mathrm{rig},L}^{\dagger,r}(V) = (S\hat{\otimes}\Bt_{\mathrm{rig},L}^{\dagger,r}) \otimes_{S\hat{\otimes}\B_K^{^\dagger,r}}D_K^{\dagger,r}(V)$.
\end{itemize}

\begin{prop}
We have 
\begin{enumerate}
\item $(\tilde{D}_L^{[r;s]}(V))^{\la} = (S\hat{\otimes}(\Bt_L^{[r;s]})^{\la}) \otimes_{S\hat{\otimes}\B_K^{^\dagger,r}}D_K^{\dagger,r}(V)$;
\item $(\tilde{D}_{\mathrm{rig},L}^{\dagger,r}(V))^{\pa} = (S\hat{\otimes}(\Bt_{\mathrm{rig},L}^{\dagger,r})^{\pa}) \otimes_{S\hat{\otimes}\B_K^{^\dagger,r}}D_K^{\dagger,r}(V)$.
\end{enumerate}
\end{prop}
\begin{proof}
The first thing to prove is that the elements of $D_K^{\dagger,r}(V)$, seen as elements of $D_K^{[r;s]}(V)$ for $s \geq r$, are locally analytic (and hence pro-analytic as elements of $D_{\mathrm{rig},K}^{\dagger,r}(V)$). By proposition \ref{prop sufficient for locana}, it suffices to check that there exists a compact open subgroup $H$ of $\Gamma_K$ such that for all $g \in H$, $||g-1|| < p^{-\frac{1}{p-1}}$ on $D_K^{[r;s]}(V)$ for $s \geq r$. By the second point of proposition \ref{prop action of gamma loc ana on phigamma}, we can take $H=\Gamma_M$.

Using this result and proposition \ref{lainla and painpa}, we get that
$$(\tilde{D}_L^{[r;s]}(V))^{\la} = (S\hat{\otimes}\Bt_L^{[r;s]})^{\la} \otimes_{S\hat{\otimes}\B_K^{^\dagger,r}}D_K^{\dagger,r}(V)$$
and that
$$(\tilde{D}_{\mathrm{rig},L}^{\dagger,r}(V))^{\pa} = (S\hat{\otimes}\Bt_{\mathrm{rig},L}^{\dagger,r})^{\pa} \otimes_{S\hat{\otimes}\B_K^{^\dagger,r}}D_K^{\dagger,r}(V).$$
We can now use proposition \ref{prop trivial action = standard loc ana}, which tells us that 
$$(S\hat{\otimes}\Bt_L^{[r;s]})^{\la}=S\hat{\otimes}(\Bt_L^{[r;s]})^{\la}$$
and that
$$(S\hat{\otimes}\Bt_{\mathrm{rig},L}^{\dagger,r})^{\pa}=S\hat{\otimes}(\Bt_{\mathrm{rig},L}^{\dagger,r})^{\pa},$$
which concludes the proof.
\end{proof}

\subsection{Monodromy descent and overconvergent families of $(\phi,\tau)$-modules}
We will now prove a monodromy descent theorem in order to produce a family of overconvergent $(\phi,\tau)$-modules attached to a family of $p$-adic representations of $\G_K$, using the overconvergent family of $(\phi,\Gamma_K)$-modules attached to it by \cite{BC08} as an input.

Let $M$ be a free $S\hat{\otimes}(\Bt_{\mathrm{rig},L}^\dagger)^{\pa}$-module of rank $d$, endowed with a surjective Frobenius $\phi : M \rightarrow M$ and with a pro-analytic action of $\Gal(L/K)$. We have:

\begin{lemm}
\label{descentI}
Let $r \geq 0$ be such that $M$ and all its structures are defined over $S\hat{\otimes}(\Bt_{\mathrm{rig},L}^{\dagger,r})^{\pa}$ and such that $b,\frac{1}{b} \in \Bt_{\mathrm{rig},L}^{\dagger,r}$. Let $m_1,\cdots,m_d$ be a basis of $M$. If $I$ is a closed interval with $I \subset [r,+\infty[$, we let $M^I = \oplus S\hat{\otimes}(\Bt_L^I)^{\la}\cdot m_i$. Then $(M^I)^{\nabla_{\gamma}=0}$ is a  $S\hat{\otimes}(\Bt_L^I)^{\nabla_{\gamma}=0}$-module free of rank $d$ such that
$$M^I = (S\hat{\otimes}\Bt_L^I)^{\la}\otimes_{(S\hat{\otimes}(\Bt_L^I)^{\la})^{\nabla_{\gamma}=0}}(M^I)^{\nabla_{\gamma}=0}.$$
\end{lemm}
\begin{proof}
Let $D_{\gamma}=\Mat(\partial_{\gamma})$. In order to prove the lemma, it suffices to show that there exists $H \in \GL_d((S\hat{\otimes}\Bt_L^I)^{\la})$ such that $\partial_{\gamma}(H)+D_{\gamma}H = 0$.

For $k \in \N$ let $D_k = \Mat(\partial_{\gamma}^k)$. For $n$ big enough, the series given by
$$H = \sum_{k \geq 0}(-1)^kD_k\frac{(b_{\gamma}-b_n^{\tau})^k}{k!}$$
converges in $M_d((S\hat{\otimes}\Bt_L^I)^{\la})$ to a solution of $\partial_{\gamma}(H)+D_{\gamma}H = 0$. Moreover, for $n$ big enough, we have $||D_k(b_{\gamma}-b_n^{\tau})^k/k!|| < 1$ for $k \geq 1$ so that $H \in \GL_d((S\hat{\otimes}\Bt_L^I)^{\la})$. 
\end{proof}

\begin{theo}
\label{descentrig}
If $M$ is a free $(S\hat{\otimes}\Bt_{\mathrm{rig},L}^{\dagger})^{\pa}$-module of rank $d$, endowed with a surjective Frobenius $\phi$ and a compatible pro-analytic action of $\Gal(L/K)$, such that $\nabla_{\gamma}(M) \subset M$, then $M^{\nabla_{\gamma} = 0}$ is a free $((S\hat{\otimes}\Bt_{\mathrm{rig},L}^{\dagger})^{\pa})^{\nabla_{\gamma}=0}$-module of rank $d$ and we have 
$$M = ((S\hat{\otimes}\Bt_{\mathrm{rig},L}^{\dagger})^{\pa})\otimes_{((S\hat{\otimes}\Bt_{\mathrm{rig},L}^{\dagger})^{\pa})^{\nabla_{\gamma}=0}}M^{\nabla_{\gamma}=0}.$$
\end{theo}
\begin{proof}
Lemma \ref{descentI} allows us to find solutions on every closed interval $I$ with $I \subset [r,+\infty[$ and we now explain how to glue these solutions using the Frobenius as in the proof of theorem 6.1 of \cite{Ber14MultiLa}.

Let $I$ be such that $I \cap pI \neq \emptyset$ and let $J = I \cap pI$. Let $m_1,\cdots,m_d$ be a basis of $(M^I)^{\nabla_{\gamma}=0}$. The Frobenius $\phi$ defines bijections $\phi^k~: (M^I)^{\nabla_{\gamma}=0} \to (M^{p^kI})^{\nabla_{\gamma}=0}$ for all $k \geq 0$. Let $P \in M_d((S\hat{\otimes}\Bt_L^J)^{\la})$ be the matrix of $(\phi(m_1),\cdots,\phi(m_d))$ in the basis $(m_1,\cdots,m_d)$.

Since $(m_1,\cdots,m_d)$ is a basis of $M^I$ by lemma \ref{descentI}, it also is a basis of $M^J$, so that $M^J = (S\hat{\otimes}\Bt_L^J)^{\la}\otimes_{S\hat{\otimes}\Bt_L^I}M^I$. But $M^I = (S\hat{\otimes}\Bt_L^I)^{\la}\otimes_{((S\hat{\otimes}\Bt_L^I)^{\la})^{\nabla_{\gamma}=0}}(M^I)^{\nabla_{\gamma}=0}$, so that
$$M^J = (S\hat{\otimes}\Bt_L^J)^{\la}\otimes_{((S\hat{\otimes}\Bt_L^I)^{\la})^{\nabla_{\gamma}=0}}(M^I)^{\nabla_{\gamma}=0}.$$
We then have
$$(M^J)^{\nabla_{\gamma}=0} = ((S\hat{\otimes}\Bt_L^J)^{\la}\otimes_{((S\hat{\otimes}\Bt_L^I)^{\la})^{\nabla_{\gamma}=0}}(M^I)^{\nabla_{\gamma}=0})^{\nabla_{\gamma}=0}$$
and thus
$$(M^J)^{\nabla_{\gamma}=0} = ((S\hat{\otimes}\Bt_L^J)^{\la})^{\nabla_{\gamma}=0}\otimes_{((S\hat{\otimes}\Bt_L^I)^{\la})^{\nabla_{\gamma}=0}}(M^I)^{\nabla_{\gamma}=0}$$
so that $(m_1,\cdots,m_d)$ is also a basis of $(M^J)^{\nabla_{\gamma}=0}$. 
For the same reasons, $(\phi(m_1),\cdots,\phi(m_d))$ is also a basis of $(M^J)^{\nabla_{\gamma}=0}$ and thus $P \in \GL_d((S\hat{\otimes}\Bt_L^I)^{\la})^{\nabla_{\gamma}=0})$.

By proposition \ref{invarnabla} we have $((S\hat{\otimes}\Bt_L^I)^{\la})^{\nabla_{\gamma}=0}) = \cup_{N,n}\B_{\tau,N,n}^I$, where $N$ runs through the finite extensions of $K$ contained in $L$. Therefore there exists a finite extension $N$ of $K$, contained in $L$, and $n \geq 0$ such that $P \in \GL_d(S\hat{\otimes}\B_{\tau,N,n}^I)$. For $k \geq 0$, let $I_k = p^kI$ and $J_k = I_k \cap I_{k+1}$, and let $E_k = \oplus_{i=1}^d\B_{\tau,N,n}^{I_k}\cdot \phi^k(m_i)$. Since $P \in \GL_d(\B_{N,n}^I)$, we have $\phi^k(P) \in \GL_d(S\hat{\otimes}\B_{\tau,N,n}^{J_k})$, and hence
$$S\hat{\otimes}\B_{\tau,N,n}^{J_k}\otimes_{S\hat{\otimes}\B_{\tau,N,n}^{I_k}}E_k = S\hat{\otimes}\B_{\tau,N,n}^{J_k}\otimes_{S\hat{\otimes}\B_{\tau,N,n}^{I_{k+1}}}E_{k+1}$$
for all $k \geq 0$. The $\{E_k\}_{k \geq 0}$ form therefore a vector bundle over $S\hat{\otimes}\B_{\tau,N,n}^{[r;+\infty[}$ for $r = \min(I)$. By theorem \ref{GlueThm} there exist $n_1,\cdots,n_d$ elements of $\cap_{k \geq 0}E_k \subset M$ such that $E_k = \oplus_{i=1}^dS\hat{\otimes}\B_{\tau,N,n}^{I_k}\cdot n_i$ for all $k \geq 0$. These elements give us a basis of $M^{\nabla_{\gamma}=0}$ over $(S\hat{\otimes}\Bt_{\mathrm{rig},N}^{\dagger})^{\pa,\nabla_{\gamma}=0}$, and thus a basis of $M$ over $(S\hat{\otimes}\Bt_{\mathrm{rig},L}^\dagger)^{\pa}$.
\end{proof}

\begin{theo}~
\label{descentriggammaexact}
Let $M$ be a free $(S\hat{\otimes}\Bt_{\mathrm{rig},L}^{\dagger})^{\pa}$-module of rank $d$, endowed with a bijective Frobenius $\phi$ and a compatible pro-analytic action of $\Gal(L/K)$, such that $\nabla_{\gamma}(M) \subset M$. Then $M^{\gamma=1}$ is a locally free $S\hat{\otimes}\B_{\tau,\mathrm{rig},K,\infty}^\dagger$-module of rank $d$ and we have 
$$M = ((\Bt_{\mathrm{rig},L}^{\dagger})^{\pa})\otimes_{\B_{\tau,\mathrm{rig},K,\infty}^\dagger}M^{\gamma=1}.$$
\end{theo}
\begin{proof}
Theorem \ref{descentrig} shows that $M^{\nabla_{\gamma}=0}$ is a free $((S\hat{\otimes}\Bt_{\mathrm{rig},L}^{\dagger})^{\pa})^{\nabla_{\gamma}=0}$-module of rank $d$, such that 
$$\Bt_{\mathrm{rig},L}^{\dagger}\otimes_{((\Bt_{\mathrm{rig},L}^{\dagger})^{\pa})^{\nabla_{\gamma}=0}}M^{\nabla_{\gamma}=0} =M$$
as $\phi$-modules over $S\hat{\otimes}\Bt_{\mathrm{rig},L}^{\dagger}$ endowed with a compatible action of $\Gal(L/K)$. By proposition \ref{invarnabla}, we have the equality $((S\hat{\otimes}\Bt_{\mathrm{rig},L}^{\dagger})^{\pa})^{\nabla_{\gamma}=0} = \bigcup_{n,N}S\hat{\otimes}\B_{\tau,\mathrm{rig},N,n}^{\dagger}$. There exists therefore a finite extension $N$ of $K$ contained in $L$, $n \geq 0$ and $s_1,\cdots,s_d$ a basis of $M^{\nabla_{\gamma}=0}$ such that $\Mat(\phi) \in \GL_d(S\hat{\otimes}\B_{\tau,\mathrm{rig},N,n}^{\dagger})$. We can always assume that $N/K$ is Galois and we do so in what follows. We let $M_N=\oplus_{i=1}^d(S\hat{\otimes}\B_{\tau,\mathrm{rig},N}^{\dagger})\cdot \phi^n(s_i)$, so that $M_N$ is a $\phi$-module over $S\hat{\otimes}\B_{\tau,\mathrm{rig},N}^{\dagger}$ such that $M^{\nabla_{\gamma}=0} = (S\hat{\otimes}\Bt_{\mathrm{rig},L}^{\dagger})^{\pa,\nabla_{\gamma}=0}\otimes_{S\hat{\otimes}\B_{\tau,\mathrm{rig},N}^{\dagger}}M_N$.

Moreover, since
$$M=(S\hat{\otimes}\Bt_{\mathrm{rig},L}^{\dagger})^\pa\otimes_{((S\hat{\otimes}\Bt_{\mathrm{rig},L}^{\dagger})^{\pa})^{\nabla_{\gamma}=0}}M^{\nabla_{\gamma}=0},$$
we get that
$$M = (S\hat{\otimes}\Bt_{\mathrm{rig},L}^{\dagger})^\pa\otimes_{S\hat{\otimes}\B_{\tau,\mathrm{rig},N}^\dagger}M_N,$$
so that we can endow $M_N$ with a structure of a $(\phi,\tau_N)$-module over $(S\hat{\otimes}\B_{\tau,\mathrm{rig}^\dagger,N},S\hat{\otimes}\Bt_{\mathrm{rig},L}^\dagger)$ endowed with an action of $\Gal(N/K)$, by defining the action of $\G_K$ on $S\hat{\otimes}\Bt_{\mathrm{rig},L}^\dagger \otimes_{S\hat{\otimes}\B_{\tau,\mathrm{rig},N}^\dagger}M_N$ as the one defined diagonally on the left handside of the tensor product
$$S\hat{\otimes}\Bt_{\mathrm{rig},L}^\dagger \otimes_{(S\hat{\otimes}\Bt_{\mathrm{rig},L}^\dagger)^\pa}M = S\hat{\otimes}\Bt_{\mathrm{rig},L}^\dagger \otimes_{S\hat{\otimes}\B_{\tau,\mathrm{rig},N}^\dagger}M_N.$$

In particular, by proposition \ref{prop classical etale descent}, $M_K:=M_N^{H_{\tau,K}}=M_N^{\gamma=1}$ is a family of $(\phi,\tau)$-modules, locally free of rank $d$ over $(S\hat{\otimes}\B_{\tau,\mathrm{rig},K},S\hat{\otimes}\Bt_{\mathrm{rig},L}^\dagger)$ such that $M = (S\hat{\otimes}\Bt_{\mathrm{rig},L}^{\dagger})^\pa\otimes_{S\hat{\otimes}\B_{\tau,\mathrm{rig},K}^\dagger}M_K$. By construction, we have $M_K \subset M^{\gamma=1}$ so that $M^{\gamma=1}$ is a family of $\phi$-module over $(S\hat{\otimes}\Bt_{\mathrm{rig},L}^\dagger)^{\pa,\gamma=1}$ of rank $d$, and thus
$$M= (S\hat{\otimes}\Bt_{\mathrm{rig},L}^\dagger)^\pa \otimes_{(S\hat{\otimes}\Bt_{\mathrm{rig},L}^\dagger)^{\pa,\gamma=1}}M^{\gamma=1}.$$
Since we have $(S\hat{\otimes}\Bt_{\mathrm{rig},L}^\dagger)^{\pa,\gamma=1} = S\hat{\otimes}\B_{\tau,\mathrm{rig},K,\infty}^\dagger$ by theorem \ref{theo loc ana basic Kummer case} and proposition \ref{prop trivial action = standard loc ana}, this implies the result.
\end{proof}

\begin{theo}
\label{theo ana families}
Let $V$ be a family of representations of $\G_K$ of rank $d$. Then there exists $s_0 \geq 0$ such that for any $s \geq s_0$, there exists a unique sub-$S\hat{\otimes}\B_{\tau,\rig,K}^{\dagger,s}$-module of $(S\hat{\otimes}\Bt_{\rig,L}^{\dagger,s})^{\Gal(L/K_\infty)}$ $D_{\tau,\rig,K}^{\dagger,s}(V)$, which is a family of $(\phi,\tau)$-modules over $(S\hat{\otimes}\B_{\tau,\rig,K}^{\dagger,s},\Bt_{\rig,L}^{\dagger,s})$ such that:
\begin{enumerate}
\item $D_{\tau,\rig,K}^{\dagger,s}(V)$ is a $S\hat{\otimes}\B_{\tau,\rig,K}^{\dagger,s}$-module locally free of rank $d$;
\item the map $(S\hat{\otimes}\Bt_{\rig}^{\dagger,s})\otimes_{S\hat{\otimes}\B_{\tau,\rig,K}^{\dagger,s}}D_{\tau,\rig,K}^{\dagger,s}(V) \rightarrow (S\hat{\otimes}\Bt_{\rig}^{\dagger,s})\otimes_S V$ is an isomorphism;
\item if $x \in \cal{X}$, the map $S/\mathfrak{m}_x\otimes_SD_{\tau,\rig,K}^{\dagger,s}(V) \rightarrow D_{\tau,\rig,K}^{\dagger,s}(V_x)$ is an isomorphism.
\end{enumerate}
\end{theo}
\begin{proof}
Let $V$ be a family of representations of $\G_K$ over $S$, of dimension $d$. Let $M = \tilde{D}_{\mathrm{rig},L}^{\dagger,s}(V)^{\pa}= (S\hat{\otimes}\Bt_{\mathrm{rig},L}^{\dagger,s})\otimes D_K^{\dagger,s}(V)$ where $D_K^{\dagger,s}(V)$ is the family of overconvergent $(\phi,\Gamma)$-modules attached to $V$ by theorem \ref{theo overconvergence phigamma}. By theorem \ref{descentriggammaexact}, $M^{\gamma=1}$ is a free $(S\hat{\otimes}\Bt_{\mathrm{rig},L}^{\dagger})^{\pa,\gamma=1}$-module of rank $d$, such that we have the following isomorphism:
$$S\hat{\otimes}\Bt_{\mathrm{rig},L}^{\dagger}\otimes_{(S\hat{\otimes}\Bt_{\mathrm{rig},L}^{\dagger})^{\pa,\gamma=1}}M^{\gamma=1} \simeq (S\hat{\otimes}\Bt_{\mathrm{rig}}^{\dagger}\otimes_{S}V)^{\G_L}$$
as families of $\phi$-modules over $S\hat{\otimes}\Bt_{\mathrm{rig},L}^{\dagger}$ endowed with a compatible action of $\Gal(L/K)$. 

By theorem \ref{theo loc ana basic Kummer case} and proposition \ref{prop trivial action = standard loc ana}, $(S\hat{\otimes}\Bt_{\mathrm{rig},L}^{\dagger})^{\pa,\gamma=1} = S\hat{\otimes}\B_{\tau,\mathrm{rig},K,\infty}^{\dagger}$. There exist therefore $n \geq 0$ and $s_1,\cdots,s_d$ a basis of $M^{\gamma=1}$ such that $\Mat(\phi) \in \GL_d(S\hat{\otimes}\B_{\tau,\mathrm{rig},K,n}^{\dagger})$. We let $D_{\tau,\mathrm{rig}}^{\dagger}=\oplus_{i=1}^d(S\hat{\otimes}\B_{\tau,\mathrm{rig},K}^{\dagger})\cdot \phi^n(s_i)$, so that $D_{\tau,\mathrm{rig}}^{\dagger}$ is a family of $\phi$-modules over $S\hat{\otimes}\B_{\tau,\mathrm{rig},K}^{\dagger}$ such that $M^{\gamma=1} = (S\hat{\otimes}\Bt_{\mathrm{rig},L}^{\dagger})^{\pa,\gamma=1}\otimes_{S\hat{\otimes}\B_{\tau,\mathrm{rig},K}^{\dagger}}D_{\tau,\mathrm{rig}}^{\dagger}$.

The module $D_{\tau,\mathrm{rig}}^\dagger$ is entirely determined by this condition: if $D_1,D_2$ are two $S\hat{\otimes}\B_{\tau,\mathrm{rig}K}^{\dagger}$-modules satisfying this condition and if $X$ is the base change matrix and $P_1,P_2$ the matrices of $\phi$, then $X \in \GL_d(S\hat{\otimes}\B_{\tau,\mathrm{rig},K,n}^{\dagger})$ for $n \gg 0$, but $X$ also satisfies $X=P_2^{-1}\phi(X)P_1$ so that $X \in \GL_d(S\hat{\otimes}\B_{\tau,\mathrm{rig},K}^{\dagger})$. 

This proves item $1$. Item $2$ follows from the isomorphism 
$$S\hat{\otimes}\Bt_{\mathrm{rig},L}^{\dagger}\otimes_{(S\hat{\otimes}\Bt_{\mathrm{rig},L}^{\dagger})^{\pa,\gamma=1}}M^{\gamma=1} \simeq (S\hat{\otimes}\Bt_{\mathrm{rig}}^{\dagger}\otimes_{S}V)^{\G_L},$$
and item $3$ follows from the unicity of the family we constructed. 
\end{proof}

\begin{rema}
Unfortunately, in contrast with the situation of \cite{BC08} and because of the method we use, we do not have any control of the $s_0$ which appears in theorem \ref{theo ana families}. 
\end{rema}

\begin{rema}
\label{rema same tau to cyclo}
The same techniques could be used to produce a family of $(\phi,\Gamma)$-modules over the cyclotomic Robba ring from a family of $(\phi,\tau)$-modules.
\end{rema}

\section{An \'{e}tale descent}

In this section, we show that the families of $(\varphi, \tau)$-modules over the Robba ring $S\hat{\otimes}\B_{\tau,\rig,K}^\dagger$ associated to a family of Galois representations descend to the bounded Robba ring $S\hat{\otimes}\B_{\tau,K}^\dagger$. This is achieved analogously to results of \cite{KL10} and \cite{Hellmann16}. 

The following is a modification of an approximation lemma due to Kedlaya and Liu (\cite[Theorem 5.2]{KL10}). (Also see \cite[Lemma 5.3]{Hellmann16} in Hellmann's work.)

\begin{lemm}
\label{lem:InvertingFamily}
Let $S$ be a Banach algebra over $\Q_p$. Let $M_S$ be a free \'{e}tale $(\varphi, \tau)$-module over $S \widehat{\otimes} \B_{\tau, K}^{\dagger}$. Suppose that there exists a basis of $M_S$ on which $\varphi - 1$ acts via a matrix whose entries have positive $p$-adic valuation. Then $$ V_S = (M_S \otimes_{S \widehat{\otimes} \B_{\tau, K}^{\dagger}} (S \widehat{\otimes}_{\Q_p} \Bt^{\dagger}))^{\varphi = 1}$$ is a free $S$-linear representation.
\end{lemm}

\begin{proof}
This follows from \cite[Theorem 5.2]{KL10} once we note two things. Namely, the statement there is written for $(\varphi, \Gamma)$-modules, but the assertion and the argument is only for the $\varphi$-action. Second thing to note is that the Frobenius in our case is a priori different, but the argument in \cite[Lemma 5.1, Theorem 5.2]{KL10} takes place over the extended Robba ring where the Frobenius matches. 
\end{proof}

The following is an analogue of \cite[Lemma 5.3]{Hellmann16}, written as a restatement of this lemma.

\begin{lemm}
\label{lem:Hellmann}
For $S$ a Banach algebra over $\Q_p$, let $\tilde{\mathcal{N}}$ be a free $\phi$-module over $S \widehat{\otimes} \Bt_{\rig}^\dagger$ of rank $d$ such that there exists a basis on which $\varphi - 1$ acts via a matrix whose entries have positive $p$-adic valuation. Then $\tilde{\mathcal{N}}^{\varphi = 1}$ is free of rank $d$ as an $S$-module.  
\end{lemm}

\begin{proof}
This is Lemma \ref{lem:InvertingFamily}.
\end{proof}

We now state our étale descent theorem:

\begin{theo}
\label{thm:EtaleDescent}

Let $\mathcal{V}$ be a family representations of $\mathcal{G}_{K}$ of rank $d$ and $D_{\tau, \mathrm{rig}}^{\dagger}(V)$ the associated family of $(\varphi, \tau)$-modules associated to it over $S\hat{\otimes}\B_{\tau,\rig,K}^\dagger$. Then there exists a model $D_{\tau,K}^{\dagger}(V)$ of  $D_{\tau, \mathrm{rig}}^{\dagger}(V)$ over $S\hat{\otimes}\B_{\tau,K}^\dagger$ such that the base extension $$ \left( D_{\tau,K}^{\dagger}(V) \otimes_{S\hat{\otimes}\B_{\tau,K}^\dagger} S\hat{\otimes}\B_{\tau,\rig,K}^\dagger \right) \rightarrow D_{\tau, \mathrm{rig}}^{\dagger}(V) $$ is an isomorphism.

\end{theo}

\begin{proof}
We argue as in \cite[Theorem 5.3]{Hellmann16}. For this purpose we briefly transport to the adic space setting. Let $X = \mathrm{Spa}(S, S^{+})$ denote the adic space corresponding to the rigid analytic space associated with $S$. For $\mathcal{N}$ a family of $(\varphi, \tau)$-modules of rank $d$ We define $$ X^{adm}_{\mathcal{N}} := \left \{ x \in X | \mathrm{dim}_{k(x)} \left( (\mathcal{N} \otimes_{S\hat{\otimes}\B_{\tau,\rig,K}^\dagger} S \widehat{\otimes}  \Bt_{\rig}^\dagger) \otimes k(x) \right)^{\varphi = 1} = d \right \}. $$

Then, we first note that for the family $D_{\tau, \mathrm{rig}}^{\dagger}(V)$, $X^{adm}_{D_{\tau, \mathrm{rig}}^{\dagger}(V)} = X$ since the family of $(\varphi, \tau)$-modules comes from the family $\mathcal{V}$ of Galois representations.

Then, we note that \cite[Lemma 7.3, Theorem 7.4]{KL10} give us, for each $x \in X$, the existence of a neighbourhood $U$ of $X$ and a local \'{e}tale descent $ D_{\tau,K}^{\dagger}(V|_{U})$ over $S\hat{\otimes}\B_{\tau,\rig,K}^\dagger$ of the family $D_{\tau, \mathrm{rig}}^{\dagger}(V|_{U})$ over $S\hat{\otimes}\B_{\tau,K}^\dagger$, on which the matrix of $\varphi - 1$ has positive $p$-adic valuation. Notice again that the assertion and argument there is written for a family of $(\varphi, \Gamma)$-modules, but only uses, and constructs a model for, the $\varphi$-action. Then, statement and proof of \cite[Theorem 5.3]{Hellmann16} shows, using Lemma \ref{lem:Hellmann}, that this gives the required descent over $X^{adm}$ (i.e. the local families glue over $X^{adm}$), which is $X$ in our case as noted. This finishes the proof.
\end{proof}

The main theorem of our paper now follows:

\begin{theo}
\label{thm phitausurconv}
Let $V$ be a family of representations of $\G_K$ of rank $d$. Then there exists $s_0 \geq 0$ such that for any $s \geq s_0$, there exists a family of $(\phi,\tau)$-modules $D_{\tau,K}^{\dagger,s}(V)$ such that:
\begin{enumerate}
\item $D_{\tau,K}^{\dagger,s}(V)$ is a $S\hat{\otimes}\B_{\tau,K}^{\dagger,s}$-module locally free of rank $d$;
\item the map $(S\hat{\otimes}\Bt^{\dagger,s})\otimes_{S\hat{\otimes}\B_{\tau,K}^{\dagger,s}}D_{\tau,K}^{\dagger,s}(V) \rightarrow (S\hat{\otimes}\Bt^{\dagger,s})\otimes_S V$ is an isomorphism;
\item if $x \in \cal{X}$, the map $S/\mathfrak{m}_x\otimes_SD_{\tau,K}^{\dagger,s}(V) \rightarrow D_{\tau,K}^{\dagger,s}(V_x)$ is an isomorphism.
\end{enumerate}
\end{theo}
\begin{proof}
Items $1$ and $2$ directly follow from theorems \ref{theo ana families} and \ref{thm:EtaleDescent}. Item $3$ follows from the unicity in theorem \ref{theo ana families}.
\end{proof}

\section{Explicit computations}
\label{sec:Explicit}
In this section, we compute some explicit families of $(\phi,\tau)$-modules in some simple cases. 
\subsection{Rank $1$ $(\phi,\tau)$-modules}
We keep the same notations as introduced in \S 1. We now assume that $K=\Qp$ and we let $K_\infty$ be a Kummer extension of $\Qp$ relative to $p$. For simplicity, we also assume that $p \neq 2$. Note that, by remark 2.1.6 of \cite{gao2016loose}, in order to completely describe the $(\phi,\tau)$-module attached to some representation $V$, it suffices to give the action of $\tau$ instead of the whole action of $\Gal(L/K)$ (this was also the original definition of $(\phi,\tau)$-modules of Caruso). 

Let $E$ be a finite extension of $\Qp$. For $\delta : \Q_p^\times \to \Q_p^\times$ a continuous character, we let $\cal{R}_E(\delta)$ denote the rank $1$ $(\phi,\Gamma)$-module over $E \otimes_{\Qp}\B_{\rig,\Qp}^\dagger$ with a basis $e_\delta$ where the actions of $\phi$ and $\Gamma$ are given by $\phi(e_\delta) = \delta(p)\cdot e_\delta$ and $\gamma(e_\delta) = \delta(\chi_{\cycl}(\gamma))\cdot e_\delta$. By \cite{colmez2010representations}, every rank $1$ $(\phi,\Gamma)$-module over $E \otimes_{\Qp}\B_{\rig,\Qp}^\dagger$ is of the form $\cal{R}_E(\delta)$ for some $\delta : \Q_p^\times \to E^\times$. 
 
Recall that we put $b = \frac{t}{\lambda} \in \At_L^+$, where $\lambda= \prod_{n \geq 0}\phi^n(\frac{[\tilde{p}]}{p}-1)$ in this setting. We have $\frac{[\tilde{p}]-p}{[\epsilon][\tilde{p}]-p} = 1 - \frac{([\epsilon]-1)[\tilde{p}]}{[\epsilon][\tilde{p}]-p}$. By \cite[2.3.3]{fontaine1994corps}, $[\epsilon][\tilde{p}]-p$ is a generator of $\ker (\theta : \Atplus \to \Atplus)$ and since $[\epsilon]-1$ is killed by theta, this implies that $\alpha:=\frac{(1-[\epsilon])}{[\epsilon][\tilde{p}]-p} \in \Atplus$. 

\begin{lemm}
\label{phitau Qp(-1)}
Let $V=\Qp(-1)$. Then the associated $(\phi,\tau)$-module admits a basis $e$ in which $\phi(e) = ([\tilde{p}]-p)\cdot e$ and $\tau(e) = \prod_{n=0}^{+\infty}\phi^n(1+\alpha [\tilde{p}])\cdot e$.
\end{lemm}
\begin{proof}
The overconvergence of $(\phi,\tau)$-modules implies in particular that the $(\phi,\tau)$-module attached to $V=\Qp(-1)$ is overconvergent and thus $(\B_\tau^\dagger\otimes_\Qp V)^{H_{\tau,\Qp}}$ is of dimension $1$ over $\B_{\tau,\Qp}^\dagger$. In particular, $(\B_\tau^\dagger\otimes_\Qp V)^{H_{\tau,K}}$ is generated by an element $z \otimes a \neq 0$, and up to dividing by an element of $\Q_p^\times$, we can assume that $a=1$. Therefore there exists $z \in \B_\tau^\dagger$, $z \neq 0$, such that for all $g \in H_{\tau,\Qp}$, $g(z) = \chi_{\cycl}(g)z$. This also implies that $\G_L$ acts trivially on $z$ so that $z \in \B_{\tau,L}^\dagger$. Let $r > 0$ be such that $z \in \B_{\tau,L}^{\dagger,r}$ and such that $1/b \in \Bt_L^{\dagger,r}$. The proof of the overconvergence of $(\phi,\tau)$-modules shows that the elements of the overconvergent $(\phi,\tau)$-module lie within $(\Bt_{\rig,L}^\dagger \otimes_{\Qp}V)^{\pa}$ and therefore $z \otimes 1$ is pro-analytic for the action of $\Gal(L/\Qp)$, and thus $z \in (\Bt_{\mathrm{rig},L}^{\dagger,r})^{\pa}$.

Now if $\gamma$ is a topological generator of $\Gal(L/K_\infty)$, we have $\gamma(b)=\chi_{\cycl}(\gamma)b$, so that $z/b \in \Bt_L^I$ is left invariant by $\gamma$. Moreover, since $z$ and $1/b$  are pro-analytic vectors of $\Bt_{\mathrm{rig},L}^{\dagger,r}$, it is still the case for $z/b$. This implies that $z/b \in (\Bt_{\mathrm{rig},L}^{\dagger,r})^{\pa,\gamma=1} \B_{\tau,\mathrm{rig},\Qp,\infty}^{\dagger,r}$ by theorem \ref{theo loc ana basic Kummer case}, so that $z/b \in \B_{\tau,\mathrm{rig},\Qp,\infty}^{\dagger,r}$. 

Thus there exists $n$ such that $z/b \in \phi^{-n}(\B_{\tau,\mathrm{rig},\Qp}^{\dagger,p^nr})$ and thus $\phi^n(z/b) \in \B_{\tau,\mathrm{rig},K}^{\dagger,p^nr}$. But $z$ and $b$ are bounded elements belonging to $\Bt_L^\dagger$ and we have $\Bt_L^\dagger \cap \B_{\tau,\mathrm{rig},\Qp}^{\dagger,p^nr} =\B_{\tau,\Qp}^{\dagger,p^nr}$, so that $\phi^n(z/b) \in \B_{\tau,\Qp}^\dagger$. Since $b= \frac{t}{\lambda}$, we have $\phi^n(t)=p^nt \in \phi^n(\lambda)\cdot \B_{\tau,L}^\dagger$, and since $\phi^n(\lambda)= \frac{1}{\prod_{k=0}^{n-1}\phi^k(E([\tilde{p}])/E(0)}\cdot \lambda$, we have $t \in \lambda \cdot \B_{\tau,L}^\dagger$.

Therefore, we have $b \in \B_{\tau,L}^\dagger$ and we can take $z=b$ as a basis of $V(-1)$. The action of $\tau$ and $\phi$ on $b$ coincide with the ones given for the basis $e$ of the $(\phi,\tau)$-module.
\end{proof}

Recall that by local class field theory, the abelianization $W_{\Qp}^{\mathrm{ab}}$ of the Weil group $W_{\Qp}$ of $\Qp$ is isomorphic to $\Q_p^\times$, so that we can see any continuous character $\delta : \Q_p^\times \to \Q_p^\times$ as a continuous character of $W_{\Qp}$. Moreover, if $\delta(p) \in \Z_p^\times$ then it extends by continuity to a character of $\G_{\Qp}$. 

Note that there is a unique way of writing $\chi_{\cycl}(g)=\omega(g)\cdot \langle \chi_{\cycl}(g) \rangle$ where $\omega(g)^{p-1}=1$ and $\langle \chi_{\cycl}(g) \rangle = 1 \mod p$. The functions are still characters of $\G_{\Qp}$ and we have the following well known result:

\begin{lemm}
Every character $\G_{\Qp} \rightarrow \Z_p^\times$ is of the form $\delta = \mu_{\beta}\cdot \omega^r \cdot \langle \chi_{\cycl} \rangle^s$ where $r \in \Z/(p-1)\Z, s \in \Zp$ and $\beta \in \Z_p^\times$. 
\end{lemm}

\begin{lemm}
\label{lemma delta(p)}
If $\delta : \Q_p^\times \rightarrow \Q_p^\times$ is trivial when restricted on $\Z_p^\times$, then the $(\phi,\tau)$-module corresponding to $\cal{R}_{\Qp}(\delta)$ admits a basis $e$ in which $\phi(e) = \delta(p)\cdot e$ and the action of $\G_L$ is trivial on $e$. 
\end{lemm}
\begin{proof}
Let $e_\delta$ be the basis of $\cal{R}_{\Qp}(\delta)$ such that $\phi(e_\delta)=\delta(p)\cdot e$ and the action of $\Gamma$ is trivial on $e_\delta$, which is the same assumption as in the lemma by local class field theory. Therefore, $e_\delta \in ((\Bt_{\rig,L}^\dagger \otimes_{\B_{\rig,\Qp}^\dagger}\cal{R}_{\Qp}(\delta))^{\pa})^{\gamma=1}$ since the action of $\G_{\Qp}$ and therefore also the action of $\G_L$ is trivial on $e_\delta$. In particular, by theorem \ref{theo loc ana basic Kummer case}, there exists $n \geq 0$ such that $\phi^n(e_\delta)$ is a basis of the $(\phi,\tau)$-module corresponding to $\cal{R}_{\Qp}(\delta)$ but then $e_\delta$ also is a basis of the $(\phi,\tau)$-module, and it satisfies the stated properties.
\end{proof}

For any $g \in \G_{\Qp}$, we have $\chi_{\cycl}(g)^{p-1} \in 1+p\Zp$. Therefore, for any $a \in \Zp$, $(\chi_{\cycl}^{p-1})^a$ has a sense as a character of $\G_{\Qp}$, and if $s = (p-1)a$ then we have $(\chi_{\cycl}^{p-1})^a = \langle \chi_{\cycl} \rangle^s$. We write $T_s$ for the $\Z_p$-adic representation of $\G_{\Qp}$ corresponding to the character $\langle \chi_{\cycl} \rangle^s$ and we let $V_s = \Qp \otimes_{\Zp}T_s$.

\begin{lemm}
\label{continuityrep}
If $s_1 = s_2 \mod p^k$ then $T_{s_1} = T_{s_2} \mod p^{k+1}$
\end{lemm}
\begin{proof}
This just follows from the fact that for any $g \in \G_{\Qp}$ we have $\langle \chi_{\cycl} \rangle^s = (\chi_{\cycl}(g)^{p-1})^{\frac{1}{p-1}s}$ and the fact that $\chi_{\cycl}(g)^{p-1} \in 1+p\Zp$.
\end{proof}

For $s \in \Zp$, we let $M_{\tau}(s)$ be the $(\phi,\tau)$-module over $(\A_{\tau,K}^\dagger,\At_L^\dagger)$ having a basis $e_s$ in which $\phi(e_s) = (1-\frac{p}{[\tilde{p}]})^s\cdot e_s$ and $\tau(e_s)=[\epsilon]^s\prod_{n=0}^{+\infty}\phi^n(1+\alpha T)^s\cdot e_s$. Note that this makes sense for $s \in \Zp$ since $[\epsilon] = (1+([\epsilon]-1))$. 

\begin{lemm}
\label{continuityphitau}
If $s_1=s_2 \mod p^k$ then $M_\tau(s_1) = M_\tau(s_2) \mod (p,[\tilde{p}])^{k+1}+([\tilde{p}])^k$. 
\end{lemm}
\begin{proof}
This follows from the fact that $(1+T)^{p^k} = 1+T^k \mod (p,T)^{k+1}$. 
\end{proof}

\begin{theo}
\label{theo rank1 phitau rep}
The $(\phi,\tau)$-module corresponding to $\delta = \mu_{\beta}\cdot \omega^r \cdot \langle \chi_{\cycl} \rangle^s$ admits a basis $e$ in which $\phi(e) = \beta \cdot [\tilde{p}]^r\cdot (1-\frac{p}{[\tilde{p}]})^{-s}\cdot e$ and $\tau(e) = [\epsilon]^{-r}\prod_{n=0}^{+\infty}\phi^n(1+\alpha [\tilde{p}])^{-s}\cdot e$.
\end{theo}
\begin{proof}
By lemma \ref{phitau Qp(-1)} and compatibility with tensor products, the $(\phi,\tau)$-module attached to $\Qp(1-p)$ admits a basis $y$ in which $\phi(y) = ([\tilde{p}]-p)^{p-1}$ and $\tau(y) = \prod_{n=0}^{+\infty}\phi^n(1+\alpha [\tilde{p}])^{p-1}\cdot y$. In the basis $z = \frac{y}{[\tilde{p}]}$, we get that $\phi(z) = (1-\frac{p}{[\tilde{p}]})^{p-1}$ and $\tau(z) = \frac{1}{[\epsilon]^{p-1}}\prod_{n=0}^{+\infty}\phi^n(1+\alpha [\tilde{p}])^{p-1}\cdot z$.

Therefore, for all $s \in \N$, the $(\phi,\tau)$-module $M_{\tau}(-s)$ is the $(\phi,\tau)$-module attached to the representation $V = \Qp((1-p)s)$. By lemmas \ref{continuityrep} and \ref{continuityphitau} and since $(p-1)\N$ is a dense subset of $\Zp$, this means that for any $s \in \Zp$, the $(\phi,\tau)$-module $M_{\tau}(-s)$ is the $(\phi,\tau)$-module over $(\A_{\tau,K}^\dagger,\At_L^\dagger)$ attached to the representation $V = \Qp((1-p)s)=T_s$. 

In particular, for $s=\frac{1}{1-p}$, the $(\phi,\tau)$-module $M_{\tau}(s)$ is the one attached to $\langle \chi_{\cycl} \rangle$. By compatibility with tensor products, lemma \ref{phitau Qp(-1)}, and by the fact that $\omega = \chi_{\cycl}\cdot\langle \chi_{\cycl} \rangle^{-1}$, we get that the $(\phi,\tau)$-module attached to $\omega$ admits a basis $e$ in which $\phi(e)= [\tilde{p}]\cdot e$ and $\tau(e) = [\epsilon]\cdot e$. 

The theorem now follows by compatibility with tensor products, lemma \ref{lemma delta(p)} and our choice of normalization of local class field theory.
\end{proof}

Theorem \ref{theo rank1 phitau rep} gives us therefore a description of every $(\phi,\tau)$-module of rank $1$.

\subsection{Trianguline $(\phi,\tau)$-modules}
In \cite{colmez2010representations}, Colmez introduced the notion of trianguline representations, which are representations whose attached $(\phi,\Gamma)$-module over the Robba ring is a successive extension of rank $1$ $(\phi,\Gamma)$-modules. Colmez then computed the $(\phi,\Gamma)$-modules attached to rank $2$ trianguline representations, and those computations played a huge part in the construction of the $p$-adic Langlands correspondence for $\GL_2(\Qp)$. 

Here, our goal is to give some description of the rank $2$ $(\phi,\tau)$-modules attached to semistable representations. As in \cite{colmez2010representations}, we can define a notion of trianguline representations, relative to the theory of $(\phi,\tau)$-modules: we say that a representation is $\tau$-trianguline if its $(\phi,\tau)$-module over $\B_{\tau,\rig,K}^\dagger$ is a successive extension of rank $1$ $(\phi,\tau)$-modules. The next proposition shows that the notion of $\tau$-trianguline representations coincides with the notion of trianguline representations of Colmez. In order to keep the notations simple, we write $D_{\tau,\rig}^\dagger(\cdot)$ the functor constructed in theorem \ref{theo ana families}, from $(\phi,\Gamma)$-modules over the Robba ring to $(\phi,\tau)$-modules over $\B_{\tau,\rig,K}^\dagger$.

\begin{prop}
\label{prop eq of tannakian cat}
A sequence $0 \longrightarrow D_1 \longrightarrow D  \longrightarrow D_2  \longrightarrow 0$ in the tannakian category of $(\phi,\Gamma)$-modules over the Robba ring is exact if and only if the sequence 
$$0 \longrightarrow D_{\tau,\rig}^\dagger(D_1) \longrightarrow D_{\tau,\rig}^\dagger(D) \longrightarrow D_{\tau,\rig}^\dagger(D_2)  \longrightarrow 0$$
is exact in the category of $(\phi,\tau)$-modules. Moreover, the first sequence is split if and only if the second one is.
\end{prop}
\begin{proof}
It suffices to prove that the functor $D_{\tau,\rig}^\dagger(\cdot)$ is an equivalence of tannakian categories. In order to do so, we introduce an other category, the one of $(\phi,\Gal(L/K))$-modules over $\Bt_{\rig,L}^\dagger$, which is the category of $\phi$-modules over $\Bt_{\rig,L}^\dagger$ endowed with a compatible continuous action of $\Gal(L/K)$. Note that the construction of our functor $D_{\tau,\rig}^\dagger(\cdot)$ is constructed first by extending the scalars of a $(\phi,\Gamma)$-module $D$ to $\Bt_{\rig,L}^\dagger$, which is therefore a $(\phi,\Gal(L/K))$-module, and then showing a to descend from the $(\phi,\Gal(L/K))$-module to $D_{\tau,\rig}^\dagger(D)$. We will also denote by $D_{\tau,\rig}^\dagger(\cdot)$ the functor obtained in \S \ref{section ana families} from the category of $(\phi,\Gal(L/K))$-modules over $\Bt_{\rig,L}^\dagger$ to the category of $(\phi,\tau)$-modules over $\B_{\tau,\rig}^\dagger$. 

It is clear from the constructions of \S \ref{section ana families} that the functor $D_{\tau,\rig}^\dagger(\cdot)$ from the category of $(\phi,\Gal(L/K))$-modules over $\Bt_{\rig,L}^\dagger$ to the category of $(\phi,\tau)$-modules over $\B_{\tau,\rig}^\dagger$ induces an equivalence of categories whose quasi-inverse is the extension of scalars to $\Bt_{\rig,L}^\dagger$. The fact that the extension of scalars to $\Bt_{\rig,L}^\dagger$ is compatible with exact sequences and tensor products implies that $D_{\tau,\rig}^\dagger(\cdot)$ is an equivalence of tannakian categories.

We could use the same proof and remark \ref{rema same tau to cyclo} in order to show that the extension of scalars to $\Bt_{\rig,L}^\dagger$ induces an equivalence of tannakian categories between the category of $(\phi,\Gamma)$-modules over the Robba ring and the category of $(\phi,\Gal(L/K))$-modules over $\Bt_{\rig,L}^\dagger$. Here, we use the fact that we can apply the Tate-Sen formalism to descend from the category of $(\phi,\Gal(L/K))$-modules over $\Bt_{\rig,L}^\dagger$ to the category of $(\phi,\Gamma)$-modules over the Robba ring, and it is clear from the constructions that the functor thus obtained is a quasi-inverse to the extension of scalars to $\Bt_{\rig,L}^\dagger$.

Therefore, the functor $D_{\tau,\rig}^\dagger(\cdot)$ is an equivalence of tannakian categories, from the category of $(\phi,\Gamma)$-modules over the Robba ring to the category of $(\phi,\tau)$-modules over $\B_{\tau,\rig}^\dagger$. 
\end{proof}

Recall that given a character $\delta: \Q_p^\times \to E^\times$, we let $\cal{R}_E(\delta)$ denote the $(\phi,\Gamma)$-module with a basis $e_\delta$ where the actions of $\phi$ and $\Gamma$ are given by $\phi(e_\delta) = \delta(p)\cdot e_\delta$ and $\gamma(e_\delta) = \delta(\chi_{\cycl}(\gamma))\cdot e_\delta$. The constructions in \S 5 to produce $(\phi,\tau)$-modules from $(\phi,\Gamma)$-modules imply that, for any character $\delta: \Q_p^\times \to E^\times$, there exists $u_{\tau,\delta} \in (\Bt_{\rig,L}^\dagger)^{\pa}$, unique mod $\B_{\tau,\Qp}^\dagger$, such that $e_{\tau,\delta}:=(u_{\tau,\delta}\otimes e_\delta)$ is a basis of the corresponding $(\phi,\tau)$-module over $\B_{\tau,\rig,K}^\dagger$. Note that the uniqueness of $u_{\tau,\delta}$ comes from the uniqueness in theorem \ref{theo ana families}, and the mod $\B_{\tau,\Qp}^\dagger$ part comes from the fact that a base change inside a rank $1$ $(\phi,\tau)$-module over $\B_{\tau,\rig,\Qp}^\dagger$ is carried on by an element of $(\B_{\tau,\rig,K}^\dagger)^\times = \B_{\tau,\Qp}^\dagger$. Note that by the same reasoning as in lemma \ref{lemma delta(p)}, we can take $u_{\tau,\delta}=1$ if $\delta_{|\Z_p^\times}=1$. We let $\cal{R}_\tau(\delta)$ denote the corresponding $(\phi,\tau)$-module over $E\otimes_{\Qp}\B_{\tau,\rig,\Qp}^\dagger$. We also let $a_{\tau,\delta} \in \B_{\tau,\rig,K}^\dagger$ and $d_{\tau,\delta} \in \Bt_{\rig,L}^\dagger$ be such that $\phi(e_{\tau,\delta})=a_{\tau,\delta}\cdot e_{\tau,\delta}$ and $\tau(e_{\tau,\delta}) = d_{\tau,\delta}\cdot e_{\tau,\delta}$.

\begin{prop}
Let $D$ be a triangular $(\phi,\Gamma)$-module over $E \otimes_{\Qp}\B_{\rig,\Qp}^\dagger$, extension of $\cal{R}_E(\delta_1)$ by $\cal{R}_E(\delta_2)$, with basis $(e_1,e_2)$ in which we have
\begin{equation*}
\Mat(\phi)= 
\begin{pmatrix}
\delta_1(p) & \alpha_D \\
0 & \delta_2(p) 
\end{pmatrix}
\end{equation*}
and 
\begin{equation*}
\Mat(\gamma)= 
\begin{pmatrix}
\delta_1(\chi_\cycl(\gamma)) & \beta_D \\
0 & \delta_2(\chi_\cycl(\gamma)) 
\end{pmatrix}.
\end{equation*}
Then there exists $c_{\tau,D} \in (\Bt_{\rig,L}^\dagger)^{\pa}$, satisfying 
$$\gamma(c_{\tau,D})\delta_1(\chi(\gamma))+\delta_2(\chi(\gamma))^{-1}u_{\tau,\delta_2}\beta_D=c_{\tau,D},$$
and such that $(u_{\tau,\delta_1}\otimes e_1, c_{\tau,D}\otimes e_1+u_{\tau,\delta_2}\otimes e_2)$ is a basis of $\D_{\rig,\tau}^\dagger(D)$, in which

\begin{equation*}
\Mat(\phi)= 
\begin{pmatrix}
a_{\tau,\delta_1} & \frac{1}{u_{\tau,\delta_1}}(\phi(c_{\tau,D})\delta_1(p)+u_{\tau,\delta_2}\delta_2(p)^{-1}a_{\tau,\delta_2}-c_{\tau,D}a_{\tau,\delta_2}) \\
0 & a_{\tau,\delta_2} 
\end{pmatrix}.
\end{equation*}
and 
\begin{equation*}
\Mat(\tau)= 
\begin{pmatrix}
d_{\tau,\delta_1} & \frac{\tau(c_{\tau,D})}{u_{\tau,\delta_1}} \\
0 & d_{\tau,\delta_2} 
\end{pmatrix}.
\end{equation*}
\end{prop}
\begin{proof}
Since $e_1(E\otimes_{\Qp}\B_{\rig,\Qp}^\dagger)$ is a sub-saturated-$(\phi,\Gamma)$-module of rank $1$ of $D$, and by construction of $u_{\tau,\delta_1}$, we have that $(u_{\tau,\delta_1}\otimes e_1)(E\otimes_{\Qp}\B_{\tau,rig,\Qp}^\dagger)$ generates a sub-saturated-$(\phi,\tau)$-module of rank $1$ of $\D_{\rig,\tau}^\dagger(D)$. Therefore, we can find $a,d \in (\Bt_{\rig,L}^\dagger)^{\pa}$, with $d$ invertible, such that $(u_{\tau,\delta_1}\otimes e_1, a\otimes e_1+d\otimes e_2)$ is a basis of the $(\phi,\tau)$-module $\D_{\rig,\tau}^\dagger(D)$. In terms of base change matrix, this implies that $
\begin{pmatrix}
u_{\tau,\delta_1} & a \\
0 & d 
\end{pmatrix}
\in \GL_2((\Bt_{\rig,L})^{\pa})$ is the base change matrix from $(e_1,e_2)$ to a basis of $\D_{\rig,\tau}^\dagger(D)$. By proposition \ref{prop eq of tannakian cat}, we know that we can choose such a basis so that the $(\phi,\tau)$-module is triangular in this basis, and can be seen as an extension of $\cal{R}_{\tau}(\delta_1)$ by $\cal{R}_{\tau}(\delta_2)$. We can therefore choose $d=u_{\tau,\delta_2}$. For $(u_{\tau,\delta_1}\otimes e_1, a\otimes e_1+d\otimes e_2)$ to be a basis of $\D_{\rig,\tau}^\dagger(D)$, we have to have
$$g(a\otimes e_1+u_{\tau,\delta_2}\otimes e_2) = a\otimes e_1+u_{\tau,\delta_2}\otimes e_2$$
for all $g \in \Gal(L/K_\infty) \simeq \Gamma$. Thus, we have
$$\gamma(a)\delta_1(\chi(\gamma))+\delta_2(\chi(\gamma))^{-1}u_{\tau,\delta_2}\beta_D=a,$$
using the fact that $\gamma(u_{\tau,\delta_2}) = \delta_2(\chi(\gamma))^{-1}u_{\tau,\delta_2}$ by definition of $u_{\tau,\delta_2}$. 

We now compute the matrices of $\phi$ and $\tau$ in the basis $(u_{\tau,\delta_1}\otimes e_1, a\otimes e_1+u_{\tau,\delta_2}\otimes e_2)$. We already know that $\phi(u_{\tau,\delta_1}\otimes e_1) = a_{\tau,\delta_1}\cdot (u_{\tau,\delta_1}\otimes e_1)$, and that $\tau(u_{\tau,\delta_1}\otimes e_1) = d_{\tau,\delta_1}\cdot (u_{\tau,\delta_1}\otimes e_1)$. We have
$$\phi(a\otimes e_1+u_{\tau,\delta_2}\otimes e_2) = \phi(a)\delta_1(p)e_1+\phi(u_{\tau,\delta_2})\otimes(\alpha_De_1+\delta_2(p)e_2)$$
and thus
$$\phi(a\otimes e_1+u_{\tau,\delta_2}\otimes e_2) = (\phi(a)\delta_1(p)+\phi(u_{\tau,\delta_2})\alpha_D)\otimes e_1+a_{\tau,\delta_2}u_{\tau,\delta_2}\otimes e_2$$
so that 
$$\phi(a\otimes e_1+u_{\tau,\delta_2}\otimes e_2)=a_{\tau,\delta_2}(a\otimes e_1+u_{\tau,\delta_2}\otimes e_2)+(\phi(a)\delta_1(p)+\phi(u_{\tau,\delta_2})\alpha_D-aa_{\tau,\delta_2})\otimes e_1.$$
Therefore, the matrix of $\phi$ in this basis is
\begin{equation*}
\Mat(\phi)= 
\begin{pmatrix}
a_{\tau,\delta_1} & \frac{1}{u_{\tau,\delta_1}}(\phi(a)\delta_1(p)+u_{\tau,\delta_2}\delta_2(p)^{-1}a_{\tau,\delta_2}-aa_{\tau,\delta_2}) \\
0 & a_{\tau,\delta_2} 
\end{pmatrix}.
\end{equation*}
For the matrix of $\tau$, we have 
$$\tau(a\otimes e_1+u_{\tau,\delta_2}\otimes e_2) = \tau(a)\otimes e_1+\tau(u_{\tau,\delta_2})\otimes e_2$$
so that
\begin{equation*}
\Mat(\tau)= 
\begin{pmatrix}
d_{\tau,\delta_1} & \frac{\tau(a)}{u_{\tau,\delta_1}} \\
0 & d_{\tau,\delta_2} 
\end{pmatrix}.
\end{equation*}
The proposition follows by taking $c_{\tau,D} := a$.
\end{proof}

Unfortunately, it is actually quite difficult to describe the action of $\Gal(L/K)$ (or even of $\tau$) for $(\phi,\tau)$-modules, because the action happens over a ring which is too big. Because of this, we want to replace the action of $\tau$ with something that acts directly on the $\phi$-module over $\B_{\tau,\rig,K}^\dagger$. We define as in \cite[\S 3]{P19bis} an operator $N_{\nabla}$ on $(\Bt_{\mathrm{rig},L}^{\dagger})^{\pa}$, by $N_{\nabla}~:= \frac{-1}{b}\nabla_{\tau}$. Since $b \in \Bt_L^{\dagger}$ and is locally analytic by \cite[Lemm. 5.1.1]{GP18}, the operator $N_{\nabla}~: (\Bt_{\mathrm{rig},L}^{\dagger})^{\pa} \rightarrow (\Bt_{\mathrm{rig},L}^{\dagger})^{\pa}$ is well defined, and more generally, the connexion $N_{\nabla}~: (\Bt_{\mathrm{rig},L}^{\dagger}\otimes D_{\tau,\mathrm{rig}}^{\dagger}(V))^{\pa} \rightarrow (\Bt_{\mathrm{rig},L}^{\dagger}\otimes D_{\tau,\mathrm{rig}}^{\dagger}(V))^{\pa}$ is well defined. Moreover, since $\nabla_{\tau}([\tilde{\pi}]) = t[\tilde{\pi}]$ and since $\lambda \in \B_{\tau,\mathrm{rig},K}^{\dagger}$, we have that $N_{\nabla}(\B_{\tau,\mathrm{rig},K}^{\dagger}) \subset \B_{\tau,\mathrm{rig},K}^{\dagger}$, and the choice of the sign is made so that the operator $N_\nabla$ we just defined on $\B_{\tau,\mathrm{rig},K}^{\dagger}$ coincides with the operator $N_{\nabla}$ defined by Kisin in \cite{KisinFiso}, because with this definition one can check that 
$N_\nabla([\tilde{\pi}])=-\lambda[\tilde{\pi}]$.

\begin{defi}
A $(\phi,N_\nabla)$-module on $\B_{\tau,\mathrm{rig},K}^{\dagger}$ is a free $\B_{\tau,\mathrm{rig},K}^{\dagger}$-module $D$ endowed with a Frobenius and a compatible operator $N_{\nabla}~: D \rightarrow D$ over $N_{\nabla}~: \B_{\tau,\mathrm{rig},K}^{\dagger} \rightarrow \B_{\tau,\mathrm{rig},K}^{\dagger}$, which means that for all $m \in D$ and for all $x \in \B_{\tau,\mathrm{rig},K}^{\dagger}$, $N_{\nabla}(x\cdot m) = N_{\nabla}(x)\cdot m +x \cdot N_{\nabla}(m)$, and wich satisfies the relation $N_\nabla \circ \phi = \frac{E([\tilde{\pi}]}{E(0)}p\phi \circ N_\nabla$.
\end{defi}

\begin{prop}
\label{stability connexion}
If $D$ is a $(\phi,\tau)$-module over $(\B_{\tau,\rig,K}^\dagger,\Bt_L^\dagger)$, then the operator $N_\nabla := \frac{-\lambda}{t}\nabla_\tau$, defined on $(\Bt_{\rig,L}^\dagger \otimes_{\B_{\tau,\rig,K}^\dagger}D)^{\pa}$ satisfies
$$N_\nabla(D) \subset D.$$ 
\end{prop}
\begin{proof}
This is \cite[Prop 3.6]{P19bis}.
\end{proof}

Given a $p$-adic representation $V$ of $\G_K$, the operator $N_\nabla$ associated with its $(\phi,\tau)$-module $D_{\tau,\rig}^\dagger(V)$ induces a structure of $(\phi,N_\nabla)$-module. Unfortunately, the functor thus obtained is no longer faithful by \cite[Prop. 3.7]{P19bis}. In the particular case of semistable representations however, one can check that the data of the $(\phi,N_\nabla)$-module is sufficient in order to recover the representation. By \cite[Prop. 4.36]{P19bis}, the $(\phi,N_\nabla)$-modules arising from $(\phi,\tau)$-modules attached to semistable representations are exactly the Breuil-Kisin modules defined in \cite{KisinFiso}. Once this identification has been made, the fact that the data of the $(\phi,N_\nabla)$-module is sufficient to recover the representation can be done through Kisin's work \cite{KisinFiso}. The following proposition gives some description of what we can expect $(\phi,N_\nabla)$-modules attached to trianguline semistable representations to look like. In what follows, we let $\lambda'$ denote $\frac{d}{d[\tilde{p}]}\lambda$. 

\begin{prop}
\label{prop what phi Nnabla look like}
Let $V$ be a trianguline semistable representation, with nonpositive Hodge-Tate weights, whose $(\phi,\Gamma)$-module is an extension of $\cal{R}(\delta_1)$ by $\cal{R}(\delta_2)$, where $\delta_1(\Z_p^\times)$ and $\delta_2(\Z_p^\times)$ belong to $\Z_p^\times$, and are respectively of weight $k_1$ and $k_2$. Then the $(\phi,N_\nabla)$-module attached to $V$ admits a basis in which
\begin{equation*}
\Mat(\phi)= 
\begin{pmatrix}
\delta_1(p)([\tilde{p}]-p)^{-k_1} & ([\tilde{p}]-p)^{\inf(-k_1,-k_2)}\alpha_V \\
0 & \delta_2(p)([\tilde{p}]-p)^{-k_2} 
\end{pmatrix}
\end{equation*}
and
\begin{equation*}
\Mat(N_\nabla)= 
\begin{pmatrix}
-k_1[\tilde{p}]\lambda' & \beta_V \\
0 & -k_2[\tilde{p}]\lambda' 
\end{pmatrix},
\end{equation*}
where $\alpha_V, \beta_V \in \B_{\tau,\rig,K}^\dagger$. Moreover, $V$ is crystalline if and only if $\beta_V = 0 \mod [\tilde{p}]$. 
\end{prop}
\begin{proof}
This is straightforward and follows directly from the fact that there exists a basis of the $(\phi,\tau)$-module attached to $V$ corresponding to the extension of the $(\phi,\tau)$-module attached to $\delta_1$ by the one attached to $\delta_2$ (which is proposition \ref{prop eq of tannakian cat}), alongside the computations of rank $1$ $(\phi,\tau)$-modules given by theorem \ref{theo rank1 phitau rep}. 

For the matrix of $N_\nabla$, we compute the operator $N_\nabla$ attached to the representation $\Qp(-k)$, with $k \geq 1$. Let $e_k$ denote a basis of $\Qp(-k)$. Then the corresponding $(\phi,\tau)$-module admits $u:=\frac{t^k}{\lambda^k}e_k$ as a basis by the same reasoning as in lemma \ref{phitau Qp(-1)}. Therefore, we have 
$$N_\nabla(u) = -\frac{\lambda}{t}\nabla_\tau(\frac{t^k}{\lambda^k})\cdot e_k = (-\frac{\lambda}{t})\cdot (-k\nabla_\tau(\lambda)\frac{t^k}{\lambda^{k+1}}\cdot e_k.$$
Thus, we get that
$$N_\nabla(u) = k\frac{1}{t}\nabla_\tau(\lambda)\cdot u,$$
which is what we wanted, because $\frac{1}{t}\nabla_\tau(\lambda)=[\tilde{p}]\frac{d}{d[\tilde{p}]}(\lambda)$. 

The rest of the proposition follows from Kisin's results \cite{KisinFiso} and once again the fact that Kisin's constructions are compatible with our definition of $(\phi,N_\nabla)$-modules thanks to \cite[Prop. 4.36]{P19bis}. Because $V$ is semistable with nonpositive weights, its corresponding $(\phi,\tau)$-module is of finite $E$-height, which implies that $\Mat(\phi) \in ([\tilde{p}]-p)^k\M_2(\B_{\tau,\rig,K}^\dagger)$. For the last condition, Kisin's theory shows that the $(\phi,N)$-module attached to semistable representations can be recovered through the $(\phi,N_\nabla)$-module, by reduction mod $[\tilde{p}]$. The operator $N$ is then the reduction mod $[\tilde{p}]$ of $N_\nabla$, and thus in our case we get that $N=0$ (which means that $V$ is crystalline) if and only if $\beta_V = 0 \mod [\tilde{p}]$. 
\end{proof}

As an example, we give a description of the $(\phi,\tau)$-module and the $(\phi,N_\nabla)$-module attached to the ``false Tate curve''. This description is a bit more explicit than the constructions above and we do not assume that $K=\Qp$ anymore. Recall that the false Tate curve $T$ can be defined as follows: it is the $\Zp$-module of rank $2$, with basis $(e_1,e_2)$ and endowed with an action of $\G_K$ given by $g(e_1)=\chi(g)e_1$ and $g(e_2)=c(g)e_1+e_2$, where $c$ is the Kummer cocycle defined in \S 1. We let $V$ be the $p$-adic representation of $\G_K$ given by $T\otimes_{\Zp}\Qp$. 

Since for all $g \in \Gal(\overline{K}/K_\infty)$, we have $c(g)=0$, this implies that $(\frac{1}{b}\cdot e_1,e_2)$ is a basis of the attached $(\phi,\tau)$-module, where $b$ is the element of \S 2 defined by $b= \frac{t}{\lambda}$. In this basis, we have
\begin{equation*}
\Mat(\phi)= 
\begin{pmatrix}
\frac{E(0)}{E([\tilde{\pi}])} & 0 \\
0 & 1
\end{pmatrix}
\end{equation*}
and
\begin{equation*}
\Mat(\tau)= 
\begin{pmatrix}
\frac{b}{\tau(b)} & b \\
0 & 1
\end{pmatrix}.
\end{equation*}  
The computations for $N_\nabla(\frac{1}{b}\cdot e_1)$ are the same as in the proof of proposition \ref{prop what phi Nnabla look like}, and for $N_\nabla(e_2)$ it suffices to note that since $\tau(e_2)=e_1+e_2$, we have $\nabla_\tau(e_2)=e_1$ and thus $N_\nabla(e_2) = -\frac{1}{b}e_1$, so that
\begin{equation*}
\Mat(N_\nabla)= 
\begin{pmatrix}
-[\tilde{\pi}]\lambda' & -1 \\
0 & 0
\end{pmatrix}.
\end{equation*} 

\bibliographystyle{amsalpha}
\bibliography{bibli}
\end{document}